\documentclass[preprint,3p,times,review]{elsarticle}

\usepackage{amssymb}
\usepackage{amsmath}
\usepackage{amsthm}

\usepackage{lineno}

\theoremstyle{remark}
\newtheorem{remark}{Remark}[section]

\usepackage[overload]{empheq}

\usepackage{subcaption}

\usepackage{xcolor}

\usepackage{booktabs}

\usepackage{hyperref}

\begin{document}

\begin{frontmatter}

%% Title, authors and addresses

%% use the tnoteref command within \title for footnotes;
%% use the tnotetext command for theassociated footnote;
%% use the fnref command within \author or \affiliation for footnotes;
%% use the fntext command for theassociated footnote;
%% use the corref command within \author for corresponding author footnotes;
%% use the cortext command for theassociated footnote;
%% use the ead command for the email address,
%% and the form \ead[url] for the home page:
%% \title{Title\tnoteref{label1}}
%% \tnotetext[label1]{}
%% \author{Name\corref{cor1}\fnref{label2}}
%% \ead{email address}
%% \ead[url]{home page}
%% \fntext[label2]{}
%% \cortext[cor1]{}
%% \affiliation{organization={},
%%             addressline={},
%%             city={},
%%             postcode={},
%%             state={},
%%             country={}}
%% \fntext[label3]{}

\title{Towards a new PGD strategy for the simulation of slender structures}

%% use optional labels to link authors explicitly to addresses:
%% \author[label1,label2]{}
%% \affiliation[label1]{organization={},
%%             addressline={},
%%             city={},
%%             postcode={},
%%             state={},
%%             country={}}
%%
%% \affiliation[label2]{organization={},
%%             addressline={},
%%             city={},
%%             postcode={},
%%             state={},
%%             country={}}

\author[Navier]{Jean Ruel\corref{cor1}} %% Author name
\ead{jean.ruel@enpc.fr}
\cortext[cor1]{Corresponding author.}

\author[Navier,INRIA]{Frédéric Legoll}
\ead{frederic.legoll@enpc.fr}

\author[Navier]{Arthur Lebée}
\ead{arthur.lebee@enpc.fr}

\author[LMPS]{Ludovic Chamoin} 
\ead{ludovic.chamoin@ens-paris-saclay.fr}

%% Author affiliation
% Jean avait mis des choses compliquées qui ne compilent pas chez moi 
\address[Navier]{Navier, ENPC, Institut Polytechnique de Paris, Université Gustave Eiffel, CNRS, Marne-la-Vallée, France}

\address[INRIA]{MATHERIALS project-team, Inria, Paris, France}

\address[LMPS]{Université Paris-Saclay, CentraleSupélec, ENS Paris-Saclay, CNRS, LMPS, Gif-sur-Yvette, France}

%% Abstract
\begin{abstract}
Effective models for slender structures derived from well-known plate (or shell) theories are justified within the limit of a small thickness, and may therefore prove limited for intermediate slenderness. On the other hand, direct 3D simulation of such structures is sub-optimal because it does not take advantage of the presence of small dimensions in some directions and is sometimes too costly and ill-conditioned. In this context, the Proper Generalized Decomposition (PGD) method, a model order reduction method based on a modal representation of the solution with separation of variables, makes it possible to obtain a 3D solution with 2D resolution complexity. In this work, an analysis of the links between the PGD reduced order model and the solution provided by plate theory is carried out using asymptotic expansion. It is shown that, in the limit of large slenderness, the first mode of the PGD exhibits Kirchhoff-Love type kinematics, but only corresponds to the asymptotic solution in very special cases of loading and boundary conditions. To capture the asymptotic solution, a new PGD strategy is introduced consisting of computing the first two modes simultaneously. We also demonstrate that the PGD is subject to shear locking, and we show how to deal with it. Numerical experiments are provided, demonstrating the interest of this approach and confirming the theoretical analysis. 
\end{abstract}

%%Graphical abstract
%%\begin{graphicalabstract}
%%\includegraphics{grabs}
%%\end{graphicalabstract}

%%Research highlights
%%\begin{highlights}
%%\item Research highlight 1
%%\item Research highlight 2
%%\end{highlights}

%% Keywords
\begin{keyword}
%% keywords here, in the form: keyword \sep keyword
Model Order Reduction \sep Proper Generalized Decomposition \sep Plate Theory \sep Asymptotic Analysis
%% PACS codes here, in the form: \PACS code \sep code

%% MSC codes here, in the form: \MSC code \sep code
%% or \MSC[2008] code \sep code (2000 is the default)

\end{keyword}

\end{frontmatter}

%% Add \usepackage{lineno} before \begin{document} and uncomment 
%% following line to enable line numbers
%\linenumbers

%% main text
%%

\section{Introduction}

Slender structures such as plates or shells are common in mechanical systems. This is particularly true in the automotive (metal sheets) and aerospace (composite panels) sectors. Numerical simulation of these three-dimensional (3D) structures, with at least one dimension smaller than the others, is therefore of current practical interest in engineering. Simplifying this problem is useful, if not necessary, to keep computational costs low. Two main strategies are available to reduce the complexity of the problem. 

The first approach is known as dimensional reduction. It involves substituting the 3D problem with a two-dimensional (2D) plate (or shell) model. This 2D model can itself be obtained using different approaches: axiomatic, asymptotic, or hierarchical. Axiomatic approaches are based on a priori kinematic or mechanical assumptions made about the 3D field, separating the out-of-plane direction from those in-plane. Asymptotic approaches are based on the explicit introduction of thickness as a small parameter going to 0 in the equations of the 3D problem. They allow to obtain plate models following a now classical procedure, and often justify axiomatic approaches in an a posteriori way \cite{ciarlet_justification_1979,ciarlet_mathematical_1997,braess_justification_2011}. In the case of elasticity, these approaches lead to well-known and widely used theories \cite{love_small_1888,reissner_effect_1945,mindlin_influence_1951,timoshenko_theory_1959}, from the simplest Kirchhoff-Love model for thin plates to more complex models valid for laminated composite plates \cite{lebee_generalization_2017,lebee_generalization_2017-1}. Despite their usefulness in many cases, these classical theories are, for the most part, justified within the limit of a small thickness and can prove limited as thickness increases. The solution provided by these models is also generally wrong near the edges of the plates, where the kinematic assumptions are not satisfied. On the other hand, in the case of highly heterogeneous structures, in the presence of nonlinearities or in the case of complex physics, it may be difficult to introduce reasonable a priori assumptions or to carry out a rigorous asymptotic analysis to reduce the dimension of the problem. We also mention hierarchical models \cite{szabo_hierarchic_1988,babuska_hierarchic_1991,babuska_hierarchic_1992}, the aim of which is to build an adaptative model to solve the 3D problem with a desired accuracy on quantities of interest. This technique involves the choice of transverse coordinate functions. Here again, in cases more complex than homogeneous or laminated plates, the optimal choice of these functions is not obvious. Simulation of the fully 3D problem is then unavoidable. Although 3D simulation is sometimes necessary, it remains sub-optimal, since it does not take advantage of the small thickness of the structure. Moreover, simulating the 3D problem may involve too many degrees of freedom and may thus request too large computational resource while being ill-conditioned.  

Model reduction techniques developed over the last few decades, which aim to reduce computational times for the simulation of complex, multi-parametric problems, are particularly attractive in this context. A variety of methods exist, differing in the way the approximation basis is constructed. They represent a second approach for reducing the complexity of the 3D plate problem. In particular, the Proper Generalized Decomposition (PGD) technique \cite{chinesta_proper_2014}, based on a modal representation of the solution with separation of variables, makes it possible to separate the plane coordinates from the out-of-plane coordinate according to the usual approach in this field. A 3D solution is thus obtained with 2D resolution complexity as shown in the seminal work \cite{bognet_advanced_2012}, taking advantage of the particular geometry of the structure. This technique, first applied to plates \cite{bognet_advanced_2012,vidal_proper_2013} and beams \cite{vidal_assessment_2012}, was next extended to shells \cite{bognet_separated_2014,vidal_shell_2014,pruliere_3d_2014}. Composite and sandwich structures were the main focus of these studies. More recently, Functional Graded Materials (FGM) composite plates have been considered \cite{kazemzadeh-parsi_proper_2021,vidal_analysis_2021}. In addition, if additional geometric or material parameters are present (Young's modulus, ply organization in a composite, etc.), they can be taken into account in the PGD decomposition \cite{vidal_explicit_2014,vidal_modeling_2017,kazemzadeh-parsi_parametric_2023,giner_proper_2013} and drastically reduce the computational costs associated with reliability analyzes \cite{gallimard_coupling_2013} or optimization problems \cite{el-ghamrawy_proper_2023}. 

More specifically, some similarities between the first mode of the PGD approximation and the solution provided by standard plate theories have been numerically observed in the literature \cite{bognet_advanced_2012}. However, to the best of the authors' knowledge, no in-depth analysis of the links between the PGD reduced model and plate theories has been provided to date.  

\medskip

In this work, an analysis of the existing connections between the PGD modes and the solution provided by standard plate theories is provided. To achieve this goal, an asymptotic analysis of the first PGD mode is performed using the formal asymptotic expansion method. The 3D linear elasticity problem is scaled before being solved by the PGD method, and the asymptotic development is constructed. It is shown that, in the limit of large slenderness, the first mode computed by the PGD approach exhibits Kirchhoff-Love kinematics but only corresponds to this model in very special cases of loading and boundary conditions. Referring to the solution given by the homogenization theory, this result can be explained by the impossibility of accurately approximating the displacement field by a single mode in the relevant energy norm. These observations motivate the following contributions of this article. 

To address the problem identified above, a modification of the standard PGD procedure is proposed, with the aim of capturing the asymptotic solution as early as the first mode computation sequence. This new PGD strategy involves computing the first two modes simultaneously. In passing, we note that the PGD is subject to shear locking, and we show how to deal with it using the selective reduced integration technique \cite{zienkiewicz_reduced_1971,hughes_simple_1977}. We also ensure that the method is able to capture boundary layers when they are present. 

Thirdly, our theoretical results are illustrated numerically on a simple but representative example of plane elasticity. The numerical results show the interest of the new PGD strategy both in terms of accuracy and computational cost, and confirm the theoretical analysis. 

\medskip

The article is organized as follows. The asymptotic development framework and the PGD reduced-order model are described in Section~\ref{sec:framework}. The asymptotic expansion procedure is detailed in Section~\ref{sec:asymptotic_expansion}, where the results obtained are also discussed. A new PGD strategy is next proposed in Section~\ref{sec:new_PGD_strategy}. Numerical experiments showing the interest of the new PGD strategy and confirming the theoretical analysis are provided in Section~\ref{sec:numerical_results}. Eventually, conclusions and perspectives are outlined in Section~\ref{sec:conclusion}. 

\section{Rescaling of the problem and PGD reduced order modeling}
\label{sec:framework}

\subsection{The 3D problem}

Unless otherwise stated, the following conventions are adopted throughout this work: Greek indices take their values in the set $\{1,2\}$, while Latin indices belong to the set $\{1,2,3\}$. It is assumed that an origin and an orthonormal basis $(\boldsymbol{e}_1,\boldsymbol{e}_2,\boldsymbol{e}_3)$ have been chosen in the 3D Euclidean space. Let $\omega^L$ be a bounded open subset of $\mathbb{R}^{2}$ of characteristic size $L$ with boundary $\partial\omega^L$ and spanned by the vectors $(\boldsymbol{e}_1,\boldsymbol{e}_2)$. 

\medskip

\begin{figure}[h!]
\centering
\includegraphics[trim=1cm 16cm 1cm 1cm, clip=true, width=\textwidth]{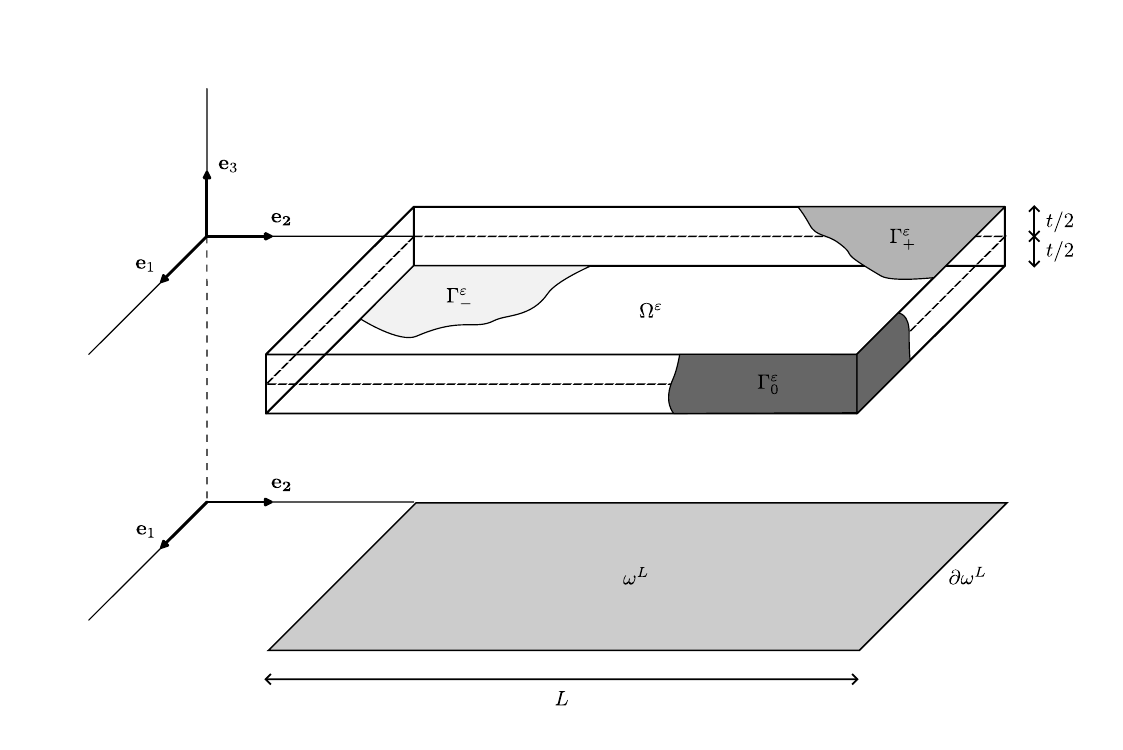}
\caption{Plate geometry}
\label{fig:plate_geo}
\end{figure}

We consider a plate of thickness $t$ occupying the domain $\Omega^{\varepsilon}$ whose points have coordinates $\boldsymbol{x}^\varepsilon=(x_1^\varepsilon,x_2^\varepsilon,x_3^\varepsilon)$ and whose boundaries $\Gamma_0^{\varepsilon}$ and $\Gamma_{\pm}^{\varepsilon}$ are defined by (see Figure~\ref{fig:plate_geo})
\begin{equation} \label{eq:plate_domain}
    \Omega^{\varepsilon} = \omega^L \times \left(-\frac{t}{2},\frac{t}{2}\right),\quad \Gamma_0^{\varepsilon} = \partial\omega^L \times \left[-\frac{t}{2},\frac{t}{2}\right],\quad \Gamma_{\pm}^{\varepsilon} = \omega^L \times \left\{\pm \frac{t}{2}\right\}.
\end{equation}
The upperscript $\displaystyle \varepsilon=\frac{t}{L}$ is the inverse of the plate slenderness, introduced at this stage to emphasize the dependence of the domain on this parameter.  

The plate is assumed to be thin ($t\ll L$), clamped on its lateral boundary $\Gamma_0^{\varepsilon}$ and subjected to body forces $\boldsymbol{f}^{\varepsilon}$ acting inside $\Omega^{\varepsilon}$ and oriented following the out-of-plane direction: 
\begin{equation*}
  \boldsymbol{f}^\varepsilon(\boldsymbol{x}^\varepsilon)~=~f_3^\varepsilon(\boldsymbol{x}^\varepsilon) \, \boldsymbol{e}_3,  
\end{equation*}
where $f_3^\varepsilon$ is assumed to be an even function of $x_3^\varepsilon$. The same density of normal surface forces
$\boldsymbol{g}^\varepsilon(x_1^\varepsilon,x_2^\varepsilon)~=~g_3^\varepsilon(x_1^\varepsilon,x_2^\varepsilon) \, \boldsymbol{e}_3$
is also applied  on the upper and lower faces $\Gamma^\varepsilon_+$ and $\Gamma^\varepsilon_-$. We thus consider a pure bending problem, but our analysis carries over to more general cases and a membrane-type problem could be investigated in a similar way. The displacement field $\boldsymbol{u}^\varepsilon= \left(u_i^\varepsilon\right)~:~\Omega^\varepsilon~\rightarrow~\mathbb{R}^3$ lives in the space
\begin{equation*} 
    V^\varepsilon = \left\{\boldsymbol{u}\in \left[H^1(\Omega^\varepsilon)\right]^3,\; \boldsymbol{u}=\boldsymbol{0}\: \text{on}\: \Gamma_0^{\varepsilon}\right\}
\end{equation*}
and is solution of the weak formulation: find $\boldsymbol{u}^\varepsilon$ in $V^\varepsilon$ such that, for any $\boldsymbol{u}^* \in V^\varepsilon$,
\begin{equation} \label{eq:epsilon_weak_form}
\begin{aligned}
    \int_{\Omega^{\varepsilon}} \boldsymbol{\epsilon}(\boldsymbol{u}^*):\boldsymbol{C}:\boldsymbol{\epsilon}(\boldsymbol{u}^\varepsilon) = \int_{\Omega^{\varepsilon}} \boldsymbol{u}^* \cdot \boldsymbol{f}^\varepsilon + \int_{\Gamma^\varepsilon_+\cup\,\Gamma^\varepsilon_-} \boldsymbol{u}^* \cdot \boldsymbol{g}^\varepsilon,
\end{aligned}
\end{equation}
where $\boldsymbol{C}$ is the fourth-order Hooke tensor and $\boldsymbol{\epsilon}(\boldsymbol{v}) = \frac{1}{2}\left(\boldsymbol{\nabla} \boldsymbol{v} + (\boldsymbol{\nabla} \boldsymbol{v})^\top \right)$ is the linearized strain tensor. This problem is equivalent to the following minimization problem, known as the principle of minimal potential energy:
\begin{equation*} 
    \boldsymbol{u} = \underset{\boldsymbol{v}\in V^\varepsilon}{\mathrm{argmin}} \ \mathcal{J}(\boldsymbol{v}),
\end{equation*}
where $\mathcal{J}$ is the potential energy defined by
\begin{equation} \label{eq:potential_energy}
    \forall\, \boldsymbol{v}\in V^\varepsilon,\quad \mathcal{J}(\boldsymbol{v}) = \frac{1}{2}\int_{\Omega^{\varepsilon}} \boldsymbol{\epsilon}(\boldsymbol{v}):\boldsymbol{C}:\boldsymbol{\epsilon}(\boldsymbol{v}) - \int_{\Omega^{\varepsilon}} \boldsymbol{v} \cdot \boldsymbol{f}^\varepsilon - \int_{\Gamma^\varepsilon_+\cup\,\Gamma^\varepsilon_-} \boldsymbol{v} \cdot \boldsymbol{g}^\varepsilon.
\end{equation}
The tensor $\boldsymbol{C}$, with coefficients $C_{ijkl}$, is symmetric, that is
\begin{equation*} \label{eq:symmetry}
     C_{ijkl} = C_{jikl} = C_{ijlk} = C_{klij},
\end{equation*}
and coercive, in the sense that there exists some constant $c_->0$ such that, for any symmetric matrix $\boldsymbol{\xi}$,
\begin{equation} \label{eq:coer}
     \boldsymbol{\xi} : \boldsymbol{C} : \boldsymbol{\xi} \geq c_- \, \boldsymbol{\xi} : \boldsymbol{\xi}.
\end{equation}
In addition, a monoclinic symmetry is assumed:
\begin{equation} \label{eq:monoclinic}
    C_{\alpha\beta\gamma 3}  = 0, \quad C_{\alpha 333}  = 0,
\end{equation}
$\boldsymbol{C}$ does not depend on $(x_1^\varepsilon,x_2^\varepsilon)$ and is an even function of $x_3^\varepsilon$. This assumption ensures the decoupling of in-plane and out-of-plane problems \cite{lebee_generalization_2017}. Note that this framework covers the case of laminated composites with a symmetrical stacking sequence. 

\medskip

Using Einstein's summation convention and the assumptions made on the tensor $\boldsymbol{C}$, we write, for any displacement fields $\boldsymbol{u}$ and $\boldsymbol{u}^*$,
\begin{equation*}
    \boldsymbol{\epsilon}^*:\boldsymbol{C}:\boldsymbol{\epsilon} = \epsilon_{\alpha\beta}^* \,C_{\alpha\beta\gamma\delta}\,\epsilon_{\gamma\delta} + \epsilon_{\alpha\beta}^*\,C_{\alpha\beta 33}\,\epsilon_{33} + 4\epsilon_{\alpha 3}^*\,C_{\alpha 3\gamma 3}\,\epsilon_{\gamma 3} + \epsilon_{33}^*\,C_{33\gamma\delta}\,\epsilon_{\gamma\delta} + \epsilon_{33}^*\,C_{3333}\,\epsilon_{33}
\end{equation*}
with
\begin{equation*}
    \epsilon_{ij} = \frac{1}{2}(u_{i,j} + u_{j,i})
\end{equation*}
and likewise for $\boldsymbol{\epsilon}^*$.

The variational problem \eqref{eq:epsilon_weak_form} therefore consists in finding $\boldsymbol{u}^\varepsilon$ in $V^\varepsilon$ such that, for any $\boldsymbol{u}^* \in V^\varepsilon$,
\begin{multline} \label{eq:weak_form_int}
  \int_{\Omega^\varepsilon} u_{\alpha,\beta}^*\, C_{\alpha\beta\gamma\delta}\,u_{\gamma,\delta}^\varepsilon + u_{\alpha,\beta}^*\, C_{\alpha\beta33}\,u_{3,3}^\varepsilon +  (u_{\alpha,3}^* + u_{3,\alpha}^*)\,C_{\alpha3\gamma3}\, (u_{\gamma,3}^\varepsilon+u_{3,\gamma}^\varepsilon) + u_{3,3}^*\, C_{33\gamma\delta}\,u_{\gamma,\delta}^\varepsilon + u_{3,3}^* \,C_{3333}\,u_{3,3}^\varepsilon \\ = \int_{\Omega^\varepsilon} u_3^*\, f_3^\varepsilon + \int_{\Gamma^\varepsilon_+\cup\,\Gamma^\varepsilon_-} u_3^* \,g_3^\varepsilon.
\end{multline}

\subsection{Scaling}

Our first objective in this work is to study the behavior of the displacement field $\boldsymbol{u}^\varepsilon$ provided by the PGD technique when the thickness of the plate decreases to 0. Since this field is defined on $\Omega^\varepsilon$, which depends itself on the thickness, we first rescale the problem so that it is posed on a fixed domain. After introducing the change of variables 
\begin{equation*}
    \boldsymbol{x} = (x_1,x_2,x_3) = \left(\frac{x_1^\varepsilon}{L}, \frac{x_2^\varepsilon}{L}, \frac{x_3^\varepsilon}{t}\right),
\end{equation*}
the plate occupies now the domain $\displaystyle \Omega = \omega \times \left(-\frac{1}{2},\frac{1}{2}\right)$ which no longer depends on $\varepsilon$. Furthermore, following the approach proposed in \cite{millet_justification_1997}, the following dimensionless fields are introduced:
\begin{equation*} 
    u_\alpha(\boldsymbol{x}) = \frac{u_\alpha^\varepsilon(\boldsymbol{x}^\varepsilon)}{U}, \quad u_3(\boldsymbol{x}) = \frac{u_3^\varepsilon(\boldsymbol{x}^\varepsilon)}{U_3}, \quad f_3(\boldsymbol{x}) = \frac{f_3^\varepsilon(\boldsymbol{x}^\varepsilon)}{F_3}, \quad g_3(x_1,x_2) = \frac{g_3^\varepsilon(x_1^\varepsilon,x_2^\varepsilon)}{G_3},\quad  A_{ijkl}(x_3) = \frac{C_{ijkl}(x_3^\varepsilon)}{E}
\end{equation*}
where the quantities $(U,U_3,F_3,G_3,E)$  are reference quantities for the problem. Since the framework of the problem is linear elasticity, small displacements are considered and $E$ is typically the Young modulus or an effective Young modulus. By defining
\begin{equation*}
    V = \left\{\boldsymbol{u}\in \left[H^1(\Omega)\right]^3,\; \boldsymbol{u}=\boldsymbol{0}\: \text{on}\: \Gamma_0\right\},
\end{equation*}
the scaling of the variational form \eqref{eq:weak_form_int} leads to find $\boldsymbol{u}$ in $V$ such that, for any $\boldsymbol{u}^* \in V$, 
\begin{equation} \label{eq:scaled_weak_form_int}
\begin{aligned}
    \upsilon^2\varepsilon^2 \int_\Omega u_{\alpha,\beta}^*\,A_{\alpha\beta\gamma\delta}\,u_{\gamma,\delta} \quad + && \upsilon\varepsilon  \int_\Omega u_{\alpha,\beta}^*\,A_{\alpha\beta33}\,u_{3,3} \quad + & \quad \upsilon^2 \int_\Omega u_{\alpha,3}^*\,A_{\alpha3\gamma3}\,u_{\gamma,3} \\
    + \quad \upsilon\varepsilon \int_\Omega u_{\alpha,3}^*\,A_{\alpha3\gamma3}\,u_{3,\gamma} \quad + && \upsilon\varepsilon \int_\Omega u_{3,\alpha}^*\,A_{\alpha3\gamma3}\,u_{\gamma,3} \quad + & \quad \varepsilon^2 \int_\Omega u_{3,\alpha}^*\,A_{\alpha3\gamma3}\,u_{3,\gamma} \\
    + \quad \upsilon\varepsilon \int_\Omega u_{3,3}^*\,A_{33\gamma\delta}\,u_{\gamma,\delta} \quad + && \int_\Omega u_{3,3}^*\,A_{3333}\,u_{3,3} \quad = & \quad \mathcal{F}_3 \int_\Omega u_3^*\, f_3 \:+\: \mathcal{G}_3 \int_{\Gamma_+\cup\,\Gamma_-} u_3^*\, g_3,
\end{aligned}
\end{equation}
where $\varepsilon$, $\upsilon$, $\mathcal{F}_3$ and $\mathcal{G}_3$ are four dimensionless numbers defined by
\begin{equation*}
    \varepsilon = \frac{t}{L}, \quad \upsilon = \frac{U}{U_3}, \quad \mathcal{F}_3 = \frac{t^2F_3}{U_3E}, \quad \mathcal{G}_3 = \frac{t\,G_3}{U_3E}.
\end{equation*}
The ratio $\varepsilon$ has already been introduced. The quotient $\upsilon$ represents the ratio between the transverse and planar displacement scales. The quotient $\mathcal{F}_3$ (resp. $\mathcal{G}_3$) can be seen as a ratio between the resultant of the body forces on the thickness (resp. the surface forces) and $E$ regarded as a reference stress. These four numbers must be linked together to obtain a problem that depends only on $\varepsilon$. The transverse displacement must remain of the order of the plate thickness, $U_3$ is therefore set equal to $t$. Linking $U$ and $U_3$ to the inclination of the mid-plane of the plate by a simple geometric argument, it is reasonable to consider $\upsilon=\varepsilon$. In order for all the applied forces to appear at the same order in Equation~\eqref{eq:scaled_weak_form_int}, we fix $\mathcal{F}_3=\mathcal{G}_3$. We resort to the choice of $\mathcal{F}_3 = \varepsilon^4$ made in \cite{millet_justification_1997}, which is such that the leading term of the displacement is of order 0 with respect to $\varepsilon$. 

\begin{remark}
This scaling is consistent with that presented in \cite{ciarlet_mathematical_1997} and justified in \cite{miara_justification_1994}. The scaling performed here can indeed be written as
\begin{equation*} 
    u_\alpha^\varepsilon(\boldsymbol{x}^\varepsilon) = \varepsilon^2 L\, u_\alpha(\boldsymbol{x}), \quad u_3^\varepsilon(\boldsymbol{x}^\varepsilon) = \varepsilon L\,u_3(\boldsymbol{x}), \quad  f_3^\varepsilon(\boldsymbol{x}^\varepsilon) = \frac{\varepsilon^3E}{L}f_3(\boldsymbol{x}), \quad  g_3^\varepsilon(x_1^\varepsilon,x_2^\varepsilon) = \varepsilon^4 E\,g_3(x_1,x_2).
\end{equation*}
\end{remark}

With the above choice of $\upsilon$, $\mathcal{F}_3$ and $\mathcal{G}_3$, the scaled problem \eqref{eq:scaled_weak_form_int} consists in finding $\boldsymbol{u}$ in $V$ such that, for any $\boldsymbol{u}^* \in V$,
\begin{multline} \label{eq:scaled_weak_form}
    \varepsilon^4  \int_\Omega u_{\alpha,\beta}^*\,A_{\alpha\beta\gamma\delta}\,u_{\gamma,\delta} \:+\: \varepsilon^2  \int_\Omega u_{\alpha,\beta}^*\, A_{\alpha\beta33}\,u_{3,3} \:+\: \varepsilon^2  \int_\Omega (u_{\alpha,3}^* \:+\: u_{3,\alpha}^*)\,A_{\alpha3\gamma3}\,(u_{\gamma,3}+u_{3,\gamma}) \:+\: \varepsilon^2  \int_\Omega u_{3,3}^*\,A_{33\gamma\delta}\,u_{\gamma,\delta} \\ + \int_\Omega u_{3,3}^*\,A_{3333}\,u_{3,3} = \varepsilon^4 \int_\Omega  u_3^*\, f_3 \:+\: \varepsilon^4 \int_{\Gamma_+\cup\,\Gamma_-}  u_3^*\, g_3.
\end{multline}

\subsection{PGD model order reduction}

We are interested here in approximating \eqref{eq:scaled_weak_form} by a PGD strategy. The principle of model order reduction using PGD is to construct a low-rank modal decomposition $\boldsymbol{u}_m$ of $\boldsymbol{u}$ in the form
\begin{equation*}
    \boldsymbol{u}_m = \sum_{k=1}^m \boldsymbol{z}_k,
\end{equation*}
where, for $1\leq k\leq m$, $\boldsymbol{z}_k$ is in the space $(V_\omega\otimes V_3)^3$ with 
\begin{equation*}
    V_\omega = H^1_0(\omega) \quad \text{and} \quad  \displaystyle V_3 = H^1\left(-\frac{1}{2},\frac{1}{2}\right).
\end{equation*}
This decomposition is classically obtained using a greedy algorithm in which each term appearing in the above sum is iteratively computed. In the asymptotic analysis that follows, we restrict ourselves to a single PGD mode. Following \cite{bognet_advanced_2012}, the displacement field $\boldsymbol{u}$ is therefore sought in the form
\begin{equation} \label{eq:PGD_decomp}
    \boldsymbol{u}(x_1,x_2,x_3) = \boldsymbol{r}(x_3)\circ \boldsymbol{v}(x_1,x_2) = 
    \begin{pmatrix}
        r_1(x_3)\, v_1(x_1,x_2) \\
        r_2(x_3)\, v_2(x_1,x_2) \\
        r_3(x_3)\, v_3(x_1,x_2)
    \end{pmatrix},
\end{equation}
where, for any $1\leq i \leq 3$, $(v_i,r_i)\in V_\omega \times V_3$ . Minimizing the potential energy functional \eqref{eq:potential_energy} on the space $(V_\omega\otimes V_3)^3$ leads to $(\boldsymbol{v},\boldsymbol{r})$ being found as solutions to the following coupled system of equations:
\begin{subequations} \label{eq:PGD_scaled_problem}
\begin{align}
     % Problème sur v
     \forall \ \boldsymbol{v}^*\in \left(V_\omega\right)^3, \quad \varepsilon^4 & \sum_{\alpha,\beta,\gamma,\delta} \int_{\omega} v_{\alpha,\beta}^* \left(\int_\mathcal{I} r_\alpha\, A_{\alpha\beta\gamma\delta}\, r_\gamma\right) v_{\gamma,\delta} \:+\: \varepsilon^2 \sum_{\alpha,\beta} \int_{\omega} v_{\alpha,\beta}^* \left(\int_\mathcal{I}  r_{\alpha}\, A_{\alpha\beta33}\, r_{3,3}\right) v_{3} \notag \\
     +\: \varepsilon^2 & \sum_{\alpha,\gamma} \int_{\omega} v_\alpha^*\left[ \left(\int_\mathcal{I} r_{\alpha,3}\, A_{\alpha 3\gamma 3}\, r_{\gamma,3}\right) v_\gamma  \:+\: \left(\int_\mathcal{I} r_{\alpha,3}\, A_{\alpha3\gamma3}\,r_{3}\right) v_{3,\gamma}\right] \notag \\
     +\: \varepsilon^2 & \sum_{\alpha,\gamma} \int_{\omega} v_{3,\alpha}^* \left[\left(\int_\mathcal{I} r_{3}\, A_{\alpha 3\gamma 3}\, r_{\gamma,3}\right) v_\gamma \:+\: \left(\int_\mathcal{I} r_3\, A_{\alpha3\gamma3}\,r_3\right)v_{3,\gamma} \right] \notag \\
     +\: \varepsilon^2 & \sum_{\gamma,\delta} \int_\omega v_{3}^* \left(\int_\mathcal{I} r_{3,3} \,A_{33\gamma\delta}\,r_{\gamma} \right)v_{\gamma,\delta} \:+\: \int_{\omega} v_3^* \left(\int_\mathcal{I} r_{3,3}\, A_{3333}\, r_{3,3}\right) v_3  \notag \\
     =\: \varepsilon^4 & \int_{\omega} v_3^* \left[\left(\int_\mathcal{I} r_3\,f_3\right) \:+\: \left(r_3^+ + r_3^-\right)g_3\right], \label{eq:PGD_scaled_problem_1} \\
    % Problème sur r
    \forall \ \boldsymbol{r}^*\in \left(V_3\right)^3, \quad \varepsilon^4 & \sum_{\alpha,\beta,\gamma,\delta} \int_\mathcal{I} r_\alpha^* \left(\int_{\omega} v_{\alpha,\beta}\, v_{\gamma,\delta}\right) A_{\alpha\beta\gamma\delta}\, r_\gamma \:+\: \varepsilon^2 \sum_{\alpha,\beta} \int_\mathcal{I} r_{\alpha}^* \left(\int_{\omega}  v_{\alpha,\beta}\,  v_{3}\right)A_{\alpha\beta33}\, r_{3,3} \notag \\
    +\: \varepsilon^2 & \sum_{\alpha,\gamma} \int_\mathcal{I} r_{\alpha,3}^*\left[\left(\int_{\omega} v_\alpha \, v_{\gamma}\right)A_{\alpha 3\gamma 3} \,r_{\gamma,3} \:+\: \left(\int_\omega v_{\alpha}\,v_{3,\gamma}\right)A_{\alpha3\gamma3}\, r_3 \right] \notag \\
    +\: \varepsilon^2 & \sum_{\alpha,\gamma} \int_\mathcal{I} r_3^* \left[\left(\int_{\omega} v_{3,\alpha}\,  v_{\gamma}\right)A_{\alpha 3\gamma 3}\, r_{\gamma,3}  \:+\: \left(\int_\omega v_{3,\alpha}\,v_{3,\gamma}\right)A_{\alpha3\gamma3}\, r_3\right] \notag \\
    +\: \varepsilon^2 & \sum_{\gamma,\delta} \int_\mathcal{I} r_{3,3}^* \left(\int_\omega v_{3}\, v_{\gamma,\delta} \right)A_{33\gamma\delta}\,r_{\gamma} \:+\: \int_\mathcal{I} r_{3,3}^* \left(\int_{\omega} v_3\,  v_3\right)A_{3333}\, r_{3,3} \notag \\ =\:\varepsilon^4 & \left[\int_\mathcal{I} r_3^* \int_{\omega} v_3\,f_3 \:+\: \left(r_3^{*+}+r_3^{*-}\right) \int_{\omega} v_3\,g_3\right], \label{eq:PGD_scaled_problem_2}
\end{align}
\end{subequations}
where $\mathcal{I}$ denotes the interval $\displaystyle \left(-\frac{1}{2},\frac{1}{2}\right)$ and $\displaystyle r_3^\pm = r_3\left(\pm\frac{1}{2}\right)$. In practice, the system of equations \eqref{eq:PGD_scaled_problem} is solved using a fixed-point algorithm. From Equation~\eqref{eq:PGD_scaled_problem}, no summation convention is used to remove any ambiguity.

\section{Asymptotic expansion procedure and discussion}
\label{sec:asymptotic_expansion}

In view of the assumption on the loading, we recall that the 3D plate problem is a pure bending problem. Consequently, the in-plane displacement $u_\alpha$ is an odd function of $x_3$ and the out-of-plane displacement $u_3$ is an even function of $x_3$. 
In terms of PGD decomposition, this implies that
\begin{equation} \label{eq:even_odd}
\text{$r_\alpha$ (resp. $r_3$) is odd (resp. even) with respect to $x_3$.}
\end{equation}

\subsection{Asymptotic expansion}

Problem~\eqref{eq:PGD_scaled_problem} depends on the inverse $\varepsilon$ of the slenderness, which is assumed to be small ($\varepsilon \ll 1$) in the case of thin plates. Our aim in this section is to formally study the limit problem arising from the application of PGD when this small parameter goes to 0. For this purpose, we assume that there exists a formal asymptotic expansion of the unknown functions:
\begin{equation} \label{eq:asymptotic_expansion}
\begin{aligned}
    \boldsymbol{v} = \boldsymbol{v}^0 + \varepsilon\, \boldsymbol{v}^1 + \varepsilon^2\, \boldsymbol{v}^2 + \dots \\
    \boldsymbol{r} = \boldsymbol{r}^0 + \varepsilon\, \boldsymbol{r}^1 + \varepsilon^2\, \boldsymbol{r}^2 + \dots
\end{aligned}
\end{equation}
where the functions $\boldsymbol{v}^p$ and $\boldsymbol{r}^p$, $p\geq 0$, are independent of $\varepsilon$. We further assume that $v_i^0$ and $r_i^0$ are non-zero for any $1\leq i\leq 3$, which can be verified a posteriori. Note that in terms of displacement, we have
\begin{equation*}
\boldsymbol{u}  = \boldsymbol{r}^0\circ \boldsymbol{v}^0 + \varepsilon\left(\boldsymbol{r}^0\circ\boldsymbol{v}^1 + \boldsymbol{r}^1\circ\boldsymbol{v}^0 \right) + \dots = \boldsymbol{u}^0 + \varepsilon\,\boldsymbol{u}^1 + \dots
\end{equation*} 
The method of formal asymptotic expansion consists in identifying the successive terms $\boldsymbol{v}^p$ and $\boldsymbol{r}^p$, $p\geq 0$, by equating the factors of the successive powers of $\varepsilon$ found in \eqref{eq:PGD_scaled_problem} when $\boldsymbol{v}$ and $\boldsymbol{r}$ are replaced by their formal expansion~\eqref{eq:asymptotic_expansion}, and by solving the resulting variational equations.

\subsubsection{Zero-order problem}

The problem of order zero is written as follows:
\begin{equation} \label{eq:0_order_problem}
\begin{aligned}
    \forall \ v_3^*\in V_\omega, \quad \int_{\omega} v_3^* \left(\int_\mathcal{I} r_{3,3}^0\, A_{3333}\, r_{3,3}^0 \right) v_3^0  & = 0, \\
    \forall \ r_3^*\in V_3, \quad \int_\mathcal{I} r_{3,3}^* \left(\int_{\omega} v_3^0\,v_3^0 \right) A_{3333}\,r_{3,3}^0 & = 0.
\end{aligned}
\end{equation}
By taking $v_3^* = v_3^0$, we deduce from the first line of \eqref{eq:0_order_problem} that
\begin{equation*}
    \int_\omega \int_\mathcal{I} A_{3333}\left(r_{3,3}^0\, v_3^0 \right)^2  = 0.
\end{equation*}
Since $A_{3333}$ is positive as a consequence of the coercivity assumption~\eqref{eq:coer}, we deduce that $r_{3,3}^0\,v_3^0=0$ which means that $u_3^0$ does not depend on $x_3$. Since $v_3^0$ is assumed to be non-zero, we have
\begin{equation} \label{eq:0_order_result}
    \boxed{r_{3,3}^0 = 0.}
\end{equation} 
The function $r_3^0$ is therefore a constant, which we denote by the same symbol $r_3^0$. The second line of \eqref{eq:0_order_problem} is also satisfied.

\subsubsection{First-order problem}

By equating the factors of order $\varepsilon$ in \eqref{eq:PGD_scaled_problem}-\eqref{eq:asymptotic_expansion}, we obtain
\begin{equation} \label{eq:1_order_problem}
\begin{aligned}
    \forall \ v_3^*\in V_\omega, \quad & \int_{\omega} v_3^* \left(\int_\mathcal{I} r_{3,3}^0\, A_{3333} \,r_{3,3}^0\right) v_3^1  \:+\:  2\int_{\omega} v_3^* \left(\int_\mathcal{I} r_{3,3}^0\, A_{3333} \,r_{3,3}^1\right) v_3^0  \:=\: 0, \\
    \forall \ r_3^*\in V_3, \quad & \int_\mathcal{I} r_{3,3}^* \left(\int_{\omega} v_3^0\,  v_3^0 \right) A_{3333}\, r_{3,3}^1  \:+\: 2\int_\mathcal{I} r_{3,3}^* \left(\int_{\omega} v_3^0\, v_3^1 \right) A_{3333}\,r_{3,3}^0  \:=\: 0.
\end{aligned}
\end{equation}
Taking into account \eqref{eq:0_order_result}, the first line of \eqref{eq:1_order_problem} is satisfied, while the second line yields, by choosing $r_3^*=r_3^1$, 
\begin{equation*} 
    \int_\mathcal{I} \int_\omega A_{3333}\left(r_{3,3}^1\,v_3^0\right)^2 = 0.
\end{equation*}
Using again the positivity of $A_{3333}$ and recalling that $v_3^0$ is non-zero by assumption, we get
\begin{equation} \label{eq:1_order_result}
    \boxed{r_{3,3}^1 = 0.}
\end{equation}
The function $r_3^1$ is thus also a constant, which we continue to denote by $r_3^1$. Combined with \eqref{eq:0_order_result}, Equation~\eqref{eq:1_order_result} shows that $u_3^1$ is also independent of $x_3$.

\subsubsection{Second-order problem}

Identifying the terms of order $\varepsilon^2$ in \eqref{eq:PGD_scaled_problem}-\eqref{eq:asymptotic_expansion}, we get 
\begin{equation*}
\begin{aligned}
     \forall \ \boldsymbol{v}^*\in \left(V_\omega\right)^3, \quad & \sum_{\alpha,\beta} \int_{\omega} v_{\alpha,\beta}^* \left(\int_\mathcal{I}  r_{\alpha}^0\, A_{\alpha\beta33}\, r_{3,3}^0\right) v_{3}^0 \:+\: \sum_{\alpha,\gamma} \int_{\omega} v_\alpha^*\left[ \left(\int_\mathcal{I} r_{\alpha,3}^0\, A_{\alpha 3\gamma 3}\, r_{\gamma,3}^0\right) v_\gamma^0  \:+\: \left(\int_\mathcal{I} r_{\alpha,3}^0\, A_{\alpha3\gamma3}\,r_{3}^0\right) v_{3,\gamma}^0\right] \\
     + &\sum_{\alpha,\gamma} \int_{\omega} v_{3,\alpha}^* \left[\left(\int_\mathcal{I} r_{3}^0\, A_{\alpha 3\gamma 3}\, r_{\gamma,3}^0\right) v_\gamma^0 \:+\: \left(\int_\mathcal{I} r_3^0 \,A_{\alpha3\gamma3}\,r_3^0\right)v_{3,\gamma}^0 \right]  \:+\:\sum_{\gamma,\delta} \int_\omega v_{3}^* \left(\int_\mathcal{I} r_{3,3}^0\, A_{33\gamma\delta}\,r_{\gamma}^0 \right)v_{\gamma,\delta}^0 \\
     + &\int_{\omega} v_3^* \left(\int_\mathcal{I} r_{3,3}^0\, A_{3333}\, r_{3,3}^0\right) v_3^2 \:+\: 2\int_{\omega} v_3^* \left(\int_\mathcal{I} r_{3,3}^0 \,A_{3333}\, r_{3,3}^2\right) v_3^0 \\
     + & \int_{\omega} v_3^* \left(\int_\mathcal{I} r_{3,3}^1\, A_{3333}\, r_{3,3}^1\right) v_3^0 \:+\: 2\int_{\omega} v_3^* \left(\int_\mathcal{I} r_{3,3}^0 \,A_{3333}\, r_{3,3}^1\right) v_3^1 \:=\: 0,
\end{aligned}
\end{equation*}
and
\begin{equation*}
\begin{aligned}
     \forall \ \boldsymbol{r}^*\in \left(V_3\right)^3, \quad & \sum_{\alpha,\beta}\int_\mathcal{I} r_{\alpha}^* \left(\int_{\omega}  v_{\alpha,\beta}^0 \,v_{3}^0\right)A_{\alpha\beta33}\, r_{3,3}^0 \:+\: \sum_{\alpha,\gamma}\int_\mathcal{I} r_{\alpha,3}^*\left[\left(\int_{\omega} v_\alpha^0\,  v_{\gamma}^0\right)A_{\alpha 3\gamma 3}\,r_{\gamma,3}^0 \:+\: \left(\int_\omega v_{\alpha}^0\,v_{3,\gamma}^0\right)A_{\alpha3\gamma3}\, r_3^0 \right] \\
     + & \sum_{\alpha,\gamma} \int_\mathcal{I} r_3^* \left[\left(\int_{\omega} v_{3,\alpha}^0\,v_{\gamma}^0\right)A_{\alpha 3\gamma 3}\,r_{\gamma,3}^0  \:+\: \left(\int_\omega v_{3,\alpha}^0\, v_{3,\gamma}^0\right)A_{\alpha3\gamma3}\,r_3^0\right] \:+\: \sum_{\gamma,\delta}\int_\mathcal{I} r_{3,3}^* \left(\int_\omega v_{3}^0\, v_{\gamma,\delta}^0 \right)A_{33\gamma\delta}\, r_{\gamma}^0 \\
     + & \int_\mathcal{I} r_{3,3}^* \left(\int_{\omega} v_3^0\,  v_3^0\right)A_{3333} \,r_{3,3}^2 \:+\: 2\int_\mathcal{I} r_{3,3}^* \left(\int_{\omega} v_3^0 \, v_3^2\right)A_{3333}\, r_{3,3}^0 \\
     + &\int_\mathcal{I} r_{3,3}^* \left(\int_{\omega} v_3^1\,  v_3^1\right)A_{3333} \,r_{3,3}^0 \:+\: 2\int_\mathcal{I} r_{3,3}^* \left(\int_{\omega} v_3^0 \, v_3^1\right)A_{3333}\, r_{3,3}^1 \:=\: 0.
\end{aligned}
\end{equation*}
Taking into account \eqref{eq:0_order_result} and \eqref{eq:1_order_result}, this problem reduces to
\begin{equation} \label{eq:2_order_problem}
\begin{aligned}
     \forall \ \boldsymbol{v}^*\in \left(V_\omega\right)^3, \quad & \sum_{\alpha,\gamma} \int_{\omega} v_\alpha^*\left[ \left(\int_\mathcal{I} r_{\alpha,3}^0\, A_{\alpha 3\gamma 3}\, r_{\gamma,3}^0\right) v_\gamma^0  \:+\: \left(\int_\mathcal{I} r_{\alpha,3}^0\, A_{\alpha3\gamma3}\,r_{3}^0\right) v_{3,\gamma}^0\right] \\
     + & \sum_{\alpha,\gamma} \int_{\omega} v_{3,\alpha}^* \left[\left(\int_\mathcal{I} r_{3}^0\, A_{\alpha 3\gamma 3}\, r_{\gamma,3}^0\right) v_\gamma^0 \:+\: \left(\int_\mathcal{I} r_3^0 \,A_{\alpha3\gamma3}\,r_3^0\right)v_{3,\gamma}^0 \right] \:=\: 0, \\
     \forall \ \boldsymbol{r}^*\in \left(V_3\right)^3, \quad & \sum_{\alpha,\gamma}\int_\mathcal{I} r_{\alpha,3}^*\left[\left(\int_{\omega} v_\alpha^0\,  v_{\gamma}^0\right)A_{\alpha 3\gamma 3}\,r_{\gamma,3}^0 \:+\: \left(\int_\omega v_{\alpha}^0\,v_{3,\gamma}^0\right)A_{\alpha3\gamma3}\, r_3^0 \right] \\
     + & \sum_{\alpha,\gamma} \int_\mathcal{I} r_3^* \left[\left(\int_{\omega} v_{3,\alpha}^0\,v_{\gamma}^0\right)A_{\alpha 3\gamma 3}\,r_{\gamma,3}^0  \:+\: \left(\int_\omega v_{3,\alpha}^0\, v_{3,\gamma}^0\right)A_{\alpha3\gamma3}\,r_3^0\right] \:+\: \sum_{\gamma,\delta}\int_\mathcal{I} r_{3,3}^* \left(\int_\omega v_{3}^0\, v_{\gamma,\delta}^0 \right)A_{33\gamma\delta}\, r_{\gamma}^0 \\
     + & \int_\mathcal{I} r_{3,3}^* \left(\int_{\omega} v_3^0\,  v_3^0\right)A_{3333} \,r_{3,3}^2 \:=\: 0.
\end{aligned}
\end{equation}
It is convenient to rewrite the first line of \eqref{eq:2_order_problem} in matrix form as follows:
\begin{equation} \label{eq:2_order_equations_1} 
    \forall \ \boldsymbol{v}^*\in \left(V_\omega\right)^3, \quad \int_\omega \int_{\mathcal{I}}\,
    \begin{pmatrix}
        r_{1,3}^0\,v_1^* + r_3^0\,v_{3,1}^* &
        r_{2,3}^0\,v_2^* + r_3^0\,v_{3,2}^*
    \end{pmatrix}\, \boldsymbol{G}\,
    \begin{pmatrix}
        r_{1,3}^0\,v_1^0 + r_3^0\,v_{3,1}^0 \\
        r_{2,3}^0\,v_2^0 + r_3^0\,v_{3,2}^0
    \end{pmatrix} = 0,
\end{equation}
where $\boldsymbol{G}$ is the out-of-plane shear stiffness matrix defined by
\begin{equation} \label{eq:shear_stiff_mat}
    \boldsymbol{G} =
    \begin{pmatrix}
        A_{1313} & A_{1323} \\
        A_{1323} & A_{2323}
    \end{pmatrix}.
\end{equation}
By denoting $\boldsymbol{d} = \begin{pmatrix} r_{1,3}^0\,v_1^0 + r_3^0\,v_{3,1}^0 & r_{2,3}^0\,v_2^0 + r_3^0\,v_{3,2}^0
\end{pmatrix}^\top$ and taking $\boldsymbol{v}^* = \boldsymbol{v}^0$ in \eqref{eq:2_order_equations_1}, we get
\begin{equation} \label{eq:2_order_int_1}
    \int_\omega \int_{\mathcal{I}}\, \boldsymbol{d}^\top \boldsymbol{G}\, \boldsymbol{d} = 0.
\end{equation}
Since the matrix $\boldsymbol{G}$ is positive definite (as a consequence of coercivity assumption~\eqref{eq:coer}), we also have
\begin{equation} \label{eq:2_order_int_2}
    \boldsymbol{d}^\top \boldsymbol{G}\, \boldsymbol{d} \geq 0.
\end{equation}
We deduce from \eqref{eq:2_order_int_1} and \eqref{eq:2_order_int_2} that
\begin{equation*} 
    \boldsymbol{d}^\top \boldsymbol{G}\, \boldsymbol{d} = 0,
\end{equation*} 
which implies that
\begin{equation} \label{eq:2_order_int_3}
    \boldsymbol{d} = 0,
\end{equation}
since $\boldsymbol{G}$ is positive definite. Equation~\eqref{eq:2_order_int_3} implies
\begin{equation} \label{eq:2_order_int_6}
    \forall\,\alpha, \quad r_{\alpha,3}^0\,v_\alpha^0 + r_3^0\, v_{3,\alpha}^0 = 0\quad \mathrm{i.e.} \quad u_{\alpha,3}^0 + u_{3,\alpha}^0 = 0.
\end{equation}
Differentiating \eqref{eq:2_order_int_6} with respect to $x_3$ and using \eqref{eq:0_order_result}, we get $r_{\alpha,33}^0\,v_\alpha^0=0$, hence $r_{\alpha,33}^0=0$ since $v_\alpha^0$ is assumed to be non-zero. We deduce that $r_\alpha^0$ is linear since $r_\alpha^0$ is odd according to~\eqref{eq:even_odd}: $r_\alpha^0(x_3) = c_\alpha \, x_3$. The constant $c_\alpha$ is non-zero (since $r_\alpha^0$ is assumed to be non-zero), and since $r_i$ and $v_i$ are defined up to a multiplicative constant, we can set for convenience $c_\alpha=r_3^0$. We thus get
\begin{equation} \label{eq:2_order_result_2}
    \boxed{r_{\alpha}^0(x_3) = r_3^0\, x_3}
\end{equation}
and, going back to \eqref{eq:2_order_int_6},
\begin{equation} \label{eq:2_order_result_1}
    \boxed{v_\alpha^0 = - v_{3,\alpha}^0.}
\end{equation}
We have thus shown that $\boldsymbol{u}^0$ has the following form:
\begin{equation} \label{eq:KL_displ}
    \boldsymbol{u}^0(x_1,x_2,x_3) = 
    \begin{pmatrix}
        -x_3\,r_3^0\,v_{3,1}^0 \\
        -x_3\,r_3^0\,v_{3,2}^0 \\
        r_3^0\,v_3^0
    \end{pmatrix} =
    \begin{pmatrix}
        -x_3\,u_{3,1}^0(x_1,x_2) \\
        -x_3\,u_{3,2}^0(x_1,x_2) \\
        u_3^0(x_1,x_2)
    \end{pmatrix},
\end{equation}
which corresponds to Kirchhoff-Love kinematics. Moreover, according to \eqref{eq:2_order_result_1}, the two functions $v_{3,\alpha}^0$ belong to the space $H^1_0(\omega)$. Since $v_3^0\in H^1_0(\omega)$, we have proved that $v_3^0$ belongs to the space $H^2_0(\omega)$. This is equivalent to writting that $v_3^0\in H^2(\omega)$ should satisfy the boundary conditions 
\begin{equation} \label{eq:bc_v}
    \left\{
        \begin{array}{rl}
            v_3^0 &= 0 \\
            \boldsymbol{\nabla} v_3^0 \cdot \boldsymbol{n} &= 0
        \end{array}
    \right. \quad \mathrm{on} \quad  \partial\omega,
\end{equation}
where $\boldsymbol{n}$ is the outer normal vector to $\partial\omega$. 
\medskip

We now turn to the second line of \eqref{eq:2_order_problem}. By successively choosing test functions $\boldsymbol{r}^*$, for which only the in-plane or out-of-plane component is non-zero, we deduce that
\begin{subequations}
\begin{align} 
    \forall \ \left(r_1^*,r_2^*\right)\in \left(V_3\right)^2, \quad   \sum_{\alpha,\gamma}& \int_\mathcal{I} r_{\alpha,3}^*  \left[\left(\int_{\omega} v_\alpha^0\,  v_{\gamma}^0\right)A_{\alpha 3\gamma 3}\,r_{\gamma,3}^0 \:+\: \left(\int_\omega v_{\alpha}^0\,v_{3,\gamma}^0\right) A_{\alpha3\gamma3}\, r_3^0 \right]  \:=\: 0,\label{eq:2_order_equations_3} \\
    \forall \ r_3^*\in V_3, \quad   \sum_{\alpha,\gamma}&\int_\mathcal{I} r_3^* \left[\left(\int_{\omega} v_{3,\alpha}^0\, v_{\gamma}^0\right)A_{\alpha 3\gamma 3}\,r_{\gamma,3}^0  \:+\: \left(\int_\omega v_{3,\alpha}^0\, v_{3,\gamma}^0\right)A_{\alpha3\gamma3}\, r_3^0 \right] \notag \\ 
    + \sum_{\gamma,\delta}& \int_\mathcal{I} r_{3,3}^* \left(\int_\omega v_{3}^0\, v_{\gamma,\delta}^0 \right)A_{33\gamma\delta}\,r_{\gamma}^0 \:+\:
     \int_\mathcal{I} r_{3,3}^* \left(\int_{\omega} v_3^0 \,v_3^0\right)A_{3333}\, r_{3,3}^2  \:=\: 0. \label{eq:2_order_equations_4} 
\end{align}
\end{subequations}
Taking into account \eqref{eq:2_order_result_1} and \eqref{eq:2_order_result_2}, we observe that Equation~\eqref{eq:2_order_equations_3} is already satisfied. Indeed, the second factor in the integrand of \eqref{eq:2_order_equations_3} reads
\begin{equation*}
    \left(\int_{\omega} v_\alpha^0\,  v_{\gamma}^0\right)A_{\alpha 3\gamma 3}\,r_{\gamma,3}^0 \:+\: \left(\int_\omega v_{\alpha}^0\,v_{3,\gamma}^0\right) A_{\alpha3\gamma3}\, r_3^0 = \left(\int_{\omega} v_\alpha^0\,  v_{\gamma}^0\right)A_{\alpha 3\gamma 3}\, r_3^0 \:-\: \left(\int_\omega v_{\alpha}^0\,v_\gamma^0\right) A_{\alpha3\gamma3}\, r_3^0 = 0.
\end{equation*}
Taking again \eqref{eq:2_order_result_1} and \eqref{eq:2_order_result_2} into account, we see that the first term in Equation~\eqref{eq:2_order_equations_4} vanishes. Indeed, the second factor in the integrand reads
\begin{equation*}
    \left(\int_{\omega} v_{3,\alpha}^0\, v_{\gamma}^0\right)A_{\alpha 3\gamma 3}\,r_{\gamma,3}^0  \:+\: \left(\int_\omega v_{3,\alpha}^0\, v_{3,\gamma}^0\right)A_{\alpha3\gamma3}\, r_3^0 = - \left(\int_{\omega} v_{3,\alpha}^0\, v_{3,\gamma}^0\right)A_{\alpha 3\gamma 3}\,r_3^0  \:+\: \left(\int_\omega v_{3,\alpha}^0\, v_{3,\gamma}^0\right)A_{\alpha3\gamma3}\, r_3^0 = 0.
\end{equation*}
Equation \eqref{eq:2_order_equations_4} thus becomes, using again \eqref{eq:2_order_result_1},
\begin{equation} \label{eq:2_order_int_7}
    \forall \ r_3^*\in V_3, \quad -\sum_{\gamma,\delta} \int_\mathcal{I} r_{3,3}^* \left( \int_\omega v_{3}^0\, v_{3,\gamma\delta}^0 \right)A_{33\gamma\delta}\,r_{\gamma}^0 \:+\: \int_\mathcal{I} r_{3,3}^* \left(\int_{\omega} \left(v_3^0\right)^2\right)A_{3333}\, r_{3,3}^2  \:=\: 0.
\end{equation}
Equation~\eqref{eq:2_order_int_7} reads in the form $\displaystyle \int_{\mathcal{I}} r_{3,3}^*\, \tau = 0$ for any $r_3^*\in V_3=H^1\left(-\frac{1}{2},\frac{1}{2}\right)$, for some function $\tau$. This implies that $\tau=0$ on $\mathcal{I}$, which reads
\begin{equation*}
    -\sum_{\gamma,\delta} \left( \int_\omega v_{3}^0\, v_{3,\gamma\delta}^0 \right)A_{33\gamma\delta}\,r_{\gamma}^0 \:+\: \left(\int_{\omega} \left(v_3^0\right)^2\right)A_{3333}\, r_{3,3}^2  \:=\: 0,
\end{equation*}
from which we obtain
\begin{equation}  \label{eq:2_order_result_3}
    \boxed{r_3^2(x_3) = \sum_{\gamma,\delta} r_3^0 \ \frac{\displaystyle \int_\omega v_{3}^0\, v_{3,\gamma\delta}^0}{\displaystyle \int_\omega \left(v_3^0\right)^2} \, \left( \int_{-1/2}^{x_3} y \, \frac{A_{33\gamma\delta}}{A_{3333}} \, \mathrm{d}y\right) \:+\: c}
\end{equation}
where $c$ is a constant. Note that the function $r_3^2$ is, up to a multiplicative factor, the Poisson effect corrector appearing in the asymptotic development of a plate (in the absence of PGD model order reduction).

\subsubsection{Third-order problem}

At this stage, our aim is to derive the equilibrium equation satisfied by $v_3^0$ which remains the only unknown function at the leading order if $r_3^0$ is fixed. The third-order problem provides some identities that will be used to establish this equilibrium equation at the next order (see Section~\ref{sec:fourth-order}). As before, we identify the terms of order $\varepsilon^3$ in \eqref{eq:PGD_scaled_problem}-\eqref{eq:asymptotic_expansion}. Taking into account \eqref{eq:0_order_result}, \eqref{eq:1_order_result}, \eqref{eq:2_order_result_1} and \eqref{eq:2_order_result_2}, we proceed as follows.

On the one hand, we consider Equation~\eqref{eq:PGD_scaled_problem_1} and observe that the first term of \eqref{eq:PGD_scaled_problem_1} does not contribute at this order, while the contribution of the second, fifth and sixth terms to the order $\varepsilon^3$ vanishes, in view of \eqref{eq:0_order_result} and \eqref{eq:1_order_result}. We are thus left with the contribution of the third and fourth terms of \eqref{eq:PGD_scaled_problem_1}, which yields 
\begin{equation} \label{eq:3_order_equations_1}
\begin{aligned}
    \forall \ \boldsymbol{v}^*\in \left(V_\omega\right)^3, \quad & \sum_{\alpha,\gamma} \int_{\omega} v_\alpha^*\left[ \left(\int_\mathcal{I} r_{\alpha,3}^0\, A_{\alpha 3\gamma 3}\, r_{\gamma,3}^0\right) v_\gamma^1 \:+\: \left(\int_\mathcal{I} r_{\alpha,3}^0\, A_{\alpha3\gamma3}\, r_{3}^0\right) v_{3,\gamma}^1\right] \\ 
    + & \sum_{\alpha,\gamma} \int_{\omega} v_\alpha^*\left[ \left(\int_\mathcal{I} r_{\alpha,3}^0\, A_{\alpha 3\gamma 3}\, r_{\gamma,3}^1\right) v_\gamma^0  \:+\: \left(\int_\mathcal{I} r_{\alpha,3}^0\, A_{\alpha3\gamma3}\,r_{3}^1\right) v_{3,\gamma}^0\right] \\ 
    + & \sum_{\alpha,\gamma} \int_{\omega} v_{3,\alpha}^* \left[\left(\int_\mathcal{I} r_{3}^0\, A_{\alpha 3\gamma3}\, r_{\gamma,3}^0\right) v_\gamma^1 \:+\: \left(\int_\mathcal{I} r_3^0\, A_{\alpha3\gamma3}\,r_3^0\right)v_{3,\gamma}^1 \right] \\ 
    + & \sum_{\alpha,\gamma} \int_{\omega} v_{3,\alpha}^* \left[\left(\int_\mathcal{I} r_{3}^0\, A_{\alpha 3\gamma 3}\, r_{\gamma,3}^1\right) v_\gamma^0 \:+\: \left(\int_\mathcal{I} r_3^0\, A_{\alpha3\gamma3}\,r_3^1\right)v_{3,\gamma}^0 \right] \:=\: 0,
\end{aligned}
\end{equation}   
where we have used the fact, in view of \eqref{eq:2_order_result_1} and \eqref{eq:2_order_result_2}, that
$$
\forall \ \left(v_1^*,v_2^*\right)\in \left(V_\omega\right)^2, \quad \sum_{\alpha,\gamma} \int_{\omega} v_\alpha^*\left[ \left(\int_\mathcal{I} r_{\alpha,3}^1\, A_{\alpha 3\gamma 3}\, r_{\gamma,3}^0\right) v_\gamma^0  \:+\: \left(\int_\mathcal{I} r_{\alpha,3}^1\, A_{\alpha3\gamma3}\,r_{3}^0\right) v_{3,\gamma}^0\right] \:=\: 0
$$
and
$$
\forall \ v_3^*\in V_\omega, \quad \sum_{\alpha,\gamma} \int_{\omega} v_{3,\alpha}^* \left[\left(\int_\mathcal{I} r_{3}^1\, A_{\alpha 3\gamma3}\, r_{\gamma,3}^0\right) v_\gamma^0 \:+\: \left(\int_\mathcal{I} r_3^1\, A_{\alpha3\gamma3}\,r_3^0\right)v_{3,\gamma}^0 \right] \:=\: 0.
$$
Using the matrix $\boldsymbol{G}$ defined by \eqref{eq:shear_stiff_mat}, we next rewrite \eqref{eq:3_order_equations_1} in matrix form as follows
\begin{equation} \label{eq:3_order_int_1}
    \forall \ \boldsymbol{v}^*\in \left(V_\omega\right)^3, \quad \int_\omega \int_{\mathcal{I}} \begin{pmatrix} r_{1,3}^0\, v_1^* + r_3^0\, v_{3,1}^* & r_{2,3}^0\, v_2^* + r_3^0\, v_{3,2}^*\end{pmatrix}\, \boldsymbol{G}\, \begin{pmatrix} r_{1,3}^0\, v_1^1 + r_{1,3}^1\, v_1^0 + r_3^0\, v_{3,1}^1 + r_3^1\, v_{3,1}^0 \\ r_{2,3}^0\, v_2^1 + r_{2,3}^1\, v_2^0 + r_3^0\, v_{3,2}^1 + r_3^1\, v_{3,2}^0 \end{pmatrix} = 0,
\end{equation}
denote $\boldsymbol{d}=\begin{pmatrix} r_{1,3}^0\, v_1^1 + r_3^0\, v_{3,1}^1 + r_3^1\, v_{3,1}^0 & r_{2,3}^0\, v_2^1 + r_3^0\, v_{3,2}^1 + r_3^1\, v_{3,2}^0 \end{pmatrix}^\top$ and $\boldsymbol{e}=\begin{pmatrix} r_{1,3}^1\, v_1^0 & r_{2,3}^1\, v_2^0 \end{pmatrix}^\top$ and take $v_\alpha^* = v_\alpha^1$ and $\displaystyle v_3^*=v_3^1 + \frac{r_3^1}{r_3^0}v_3^0$ in \eqref{eq:3_order_int_1}. This yields 
\begin{equation} \label{eq:3_order_int_2}
    \int_\omega \int_{\mathcal{I}} \boldsymbol{d}^\top \, \boldsymbol{G}\, (\boldsymbol{d}+\boldsymbol{e}) = 0.
\end{equation}

On the other hand, we consider Equation~\eqref{eq:PGD_scaled_problem_2} with a test function $\boldsymbol{r}^*$ such that only its in-plane component does not vanish, and observe that the first term of \eqref{eq:PGD_scaled_problem_2} does not contribute at this order, while the contribution of the second term to the order $\varepsilon^3$ vanishes, in view of \eqref{eq:0_order_result} and \eqref{eq:1_order_result}. We are thus left with the contribution of the third term of \eqref{eq:PGD_scaled_problem_2}, which yields 
\begin{equation} \label{eq:3_order_equations_2}
\begin{aligned}
    \forall \ \left(r_1^*,r_2^*\right)\in \left(V_3\right)^2, \quad & \sum_{\alpha,\gamma} \int_\mathcal{I} r_{\alpha,3}^*\left[\left(\int_{\omega} v_\alpha^0 \,v_{\gamma}^0\right)A_{\alpha 3\gamma 3}\,r_{\gamma,3}^1 \:+\: \left(\int_\omega v_{\alpha}^0\, v_{3,\gamma}^0\right)A_{\alpha3\gamma3}\, r_3^1 \right] \\ 
    + & \sum_{\alpha,\gamma} \int_\mathcal{I} r_{\alpha,3}^*\left[\left(\int_{\omega} v_\alpha^0\, v_{\gamma}^1\right)A_{\alpha 3\gamma 3}\, r_{\gamma,3}^0 \:+\: \left(\int_\omega v_{\alpha}^0\,v_{3,\gamma}^1\right) A_{\alpha3\gamma 3}\, r_3^0 \right] \:=\: 0,
\end{aligned}
\end{equation}
where we have used the fact that
$$
\forall \ \left(r_1^*,r_2^*\right)\in \left(V_3\right)^2, \quad \sum_\gamma \int_\mathcal{I} r_{\alpha,3}^*\left[\left(\int_{\omega} v_\alpha^1\, v_{\gamma}^0\right) A_{\alpha 3\gamma 3}\, r_{\gamma,3}^0 \:+\: \left(\int_\omega v_{\alpha}^1 \,v_{3,\gamma}^0\right) A_{\alpha3\gamma 3}\, r_3^0 \right] \:=\: 0,
$$
which is a consequence of \eqref{eq:2_order_result_1} and \eqref{eq:2_order_result_2}. We can write \eqref{eq:3_order_equations_2} in matrix form as
\begin{equation} \label{eq:3_order_int_3}
    \forall \ \left(r_1^*,r_2^*\right)\in \left(V_3\right)^2, \quad \int_{\mathcal{I}} \int_\omega \begin{pmatrix} r_{1,3}^*\, v_1^0 & r_{2,3}^*\, v_2^0 \end{pmatrix} \boldsymbol{G} \begin{pmatrix} r_{1,3}^0\, v_1^1 + r_{1,3}^1\, v_1^0 + r_3^0\, v_{3,1}^1 + r_3^1\, v_{3,1}^0 \\ r_{2,3}^0\, v_2^1 + r_{2,3}^1\, v_2^0 + r_3^0\, v_{3,2}^1 + r_3^1\, v_{3,2}^0 \end{pmatrix} = 0,
\end{equation}
with $\boldsymbol{G}$ again defined by \eqref{eq:shear_stiff_mat}. Taking $r_\alpha^*=r_\alpha^1$, Equation~\eqref{eq:3_order_int_3} gives
\begin{equation} \label{eq:3_order_int_4}
    \int_{\mathcal{I}} \int_\omega \boldsymbol{e}^\top\,\boldsymbol{G}\,(\boldsymbol{d}+\boldsymbol{e}) = 0,
\end{equation}
where $\boldsymbol{d}$ and $\boldsymbol{e}$ are defined below \eqref{eq:3_order_int_1}.
Summing \eqref{eq:3_order_int_2} and \eqref{eq:3_order_int_4}, we get
\begin{equation} \label{eq:3_order_int_5}
    \int_\omega \int_{\mathcal{I}}  (\boldsymbol{d}+\boldsymbol{e})^\top\,\boldsymbol{G}\,(\boldsymbol{d}+\boldsymbol{e}) = 0.
\end{equation}
Following the same arguments as those to deduce \eqref{eq:2_order_int_3} from \eqref{eq:2_order_int_1}, we deduce from \eqref{eq:3_order_int_5} that $\boldsymbol{d}+\boldsymbol{e}=0$, which reads
\begin{equation*}
    \forall \ \alpha, \quad r_{\alpha,3}^0\, v_\alpha^1 + r_{\alpha,3}^1\, v_\alpha^0 + r_3^0\, v_{3,\alpha}^1 + r_3^1\, v_{3,\alpha}^0 = 0 \quad \mathrm{i.e.} \quad u_{\alpha,3}^1 + u_{3,\alpha}^1 = 0,
\end{equation*}
that is, using \eqref{eq:2_order_result_1} and \eqref{eq:2_order_result_2},
\begin{equation} \label{eq:3_order_result}
\boxed{v_\alpha^1 + v_{3,\alpha}^1 = \frac{r_{\alpha,3}^1 - r_3^1}{r_3^0}v_{3,\alpha}^0.}
\end{equation}
Since $r_i$ and $v_i$ are functions with separate variables, it also follows from \eqref{eq:3_order_result}, \eqref{eq:0_order_result} and \eqref{eq:1_order_result} that \fbox{$r_{\alpha,3}^1$ is a constant function.} 

\begin{remark}
    We can also consider Equation~\eqref{eq:PGD_scaled_problem_2} with a test function such that only its out-of-plane component is non-zero. However, the information obtained is not useful for the rest of this analysis. 
\end{remark}

\subsubsection{Fourth-order problem} \label{sec:fourth-order}

To obtain the equilibrium equation on $v_3^0$, we identify the terms of order $\varepsilon^4$ in \eqref{eq:PGD_scaled_problem}-\eqref{eq:asymptotic_expansion} and successively choose test functions $\boldsymbol{v}^*$, only one component of which is non-zero. 

Taking into account \eqref{eq:0_order_result} and \eqref{eq:1_order_result}, considering \eqref{eq:PGD_scaled_problem_1} and a test function $\boldsymbol{v}^*$ such that only its in-plane component does not vanish, we get
\begin{equation} \label{eq:4_order_int_1}
\begin{aligned}
    \forall \ \left(v_1^*,v_2^*\right)\in \left(V_\omega\right)^2, \quad & \sum_{\alpha,\beta,\gamma,\delta} \int_{\omega} v_{\alpha,\beta}^* \left(\int_\mathcal{I} r_{\alpha}^0\, A_{\alpha\beta\gamma\delta}\, r_{\gamma}^0 \right) v_{\gamma,\delta}^0  \:+\: \sum_{\alpha,\beta} \int_\omega v_{\alpha,\beta}^* \left(\int_\mathcal{I} r_\alpha^0\, A_{\alpha\beta33}\,r_{3,3}^2\right)v_3^0  \\
    + & \sum_{\alpha,\gamma} \int_{\omega} v_\alpha^*\left[ \left(\int_\mathcal{I} r_{\alpha,3}^0\, A_{\alpha 3\gamma 3}\, r_{\gamma,3}^0\right) v_\gamma^2  \:+\: \left(\int_\mathcal{I} r_{\alpha,3}^0\, A_{\alpha3\gamma3}\,r_{3}^0\right) v_{3,\gamma}^2\right]  \\
    + & \sum_{\alpha,\gamma} \int_{\omega} v_\alpha^*\left[ \left(\int_\mathcal{I} r_{\alpha,3}^0\, A_{\alpha 3\gamma 3}\, r_{\gamma,3}^2\right) v_\gamma^0  \:+\: \left(\int_\mathcal{I} r_{\alpha,3}^0\, A_{\alpha3\gamma3}\,r_{3}^2\right) v_{3,\gamma}^0\right]  \\
    + & \sum_{\alpha,\gamma} \int_{\omega} v_\alpha^*\left[ \left(\int_\mathcal{I} r_{\alpha,3}^2\, A_{\alpha 3\gamma 3}\, r_{\gamma,3}^0\right) v_\gamma^0  \:+\: \left(\int_\mathcal{I} r_{\alpha,3}^2\, A_{\alpha3\gamma3}\,r_{3}^0\right) v_{3,\gamma}^0\right]  \\
    + & \sum_{\alpha,\gamma} \int_{\omega} v_\alpha^*\left[ \left(\int_\mathcal{I} r_{\alpha,3}^0\, A_{\alpha 3\gamma 3}\, r_{\gamma,3}^1\right) v_\gamma^1  \:+\: \left(\int_\mathcal{I} r_{\alpha,3}^0\, A_{\alpha3\gamma3}\,r_{3}^1\right) v_{3,\gamma}^1\right]  \\
    + & \sum_{\alpha,\gamma} \int_{\omega} v_\alpha^*\left[ \left(\int_\mathcal{I} r_{\alpha,3}^1\, A_{\alpha 3\gamma 3}\, r_{\gamma,3}^0\right) v_\gamma^1  \:+\: \left(\int_\mathcal{I} r_{\alpha,3}^1\, A_{\alpha3\gamma3}\,r_{3}^0\right) v_{3,\gamma}^1\right]  \\
    + & \sum_{\alpha,\gamma} \int_{\omega} v_\alpha^*\left[ \left(\int_\mathcal{I} r_{\alpha,3}^1\, A_{\alpha 3\gamma 3}\, r_{\gamma,3}^1\right) v_\gamma^0  \:+\: \left(\int_\mathcal{I} r_{\alpha,3}^1\, A_{\alpha3\gamma3}\,r_{3}^1\right) v_{3,\gamma}^0\right] \: =\: 0.
\end{aligned}
\end{equation}
The fifth sum of \eqref{eq:4_order_int_1} vanishes in view of \eqref{eq:2_order_result_1} and \eqref{eq:2_order_result_2}, while the sum of the last two terms is written as
\begin{equation} \label{eq:4_order_int_2}
    \sum_{\alpha,\gamma} \int_{\omega} v_\alpha^*\left[ \left(\int_\mathcal{I} r_{\alpha,3}^1\, A_{\alpha 3\gamma 3}\, r_3^0\right) \left(v_\gamma^1+ v_{3,\gamma}^1\right) \:+\: \left(\left(\int_\mathcal{I} r_{\alpha,3}^1\, A_{\alpha3\gamma3}\,r_{3}^1\right)-\left(\int_\mathcal{I} r_{\alpha,3}^1\, A_{\alpha 3\gamma 3}\, r_{\gamma,3}^1\right)\right) v_{3,\gamma}^0\right].
\end{equation}
The contribution~\eqref{eq:4_order_int_2} vanishes since its second factor in the integrand reads, after taking into account \eqref{eq:3_order_result},
\begin{multline*}
    \left(\int_\mathcal{I} r_{\alpha,3}^1\, A_{\alpha 3\gamma 3}\, r_3^0\right) \left(v_\gamma^1+ v_{3,\gamma}^1\right) \:+\: \left(\left(\int_\mathcal{I} r_{\alpha,3}^1\, A_{\alpha3\gamma3}\,r_{3}^1\right)-\left(\int_\mathcal{I} r_{\alpha,3}^1\, A_{\alpha 3\gamma 3}\, r_{\gamma,3}^1\right)\right) v_{3,\gamma}^0 \\= \left(\int_\mathcal{I} A_{\alpha 3\gamma 3}\right) r_{\alpha,3}^1\left(r_{\gamma,3}^1-r_3^1\right)v_{3,\gamma}^0 \:+\: \left(\int_\mathcal{I} A_{\alpha 3\gamma 3}\right) r_{\alpha,3}^1\left(r_3^1-r_{\gamma,3}^1\right)v_{3,\gamma}^0=0.
\end{multline*}
Thus, Equation~\eqref{eq:4_order_int_1} becomes, after by replacing $r_{\alpha,3}^0$ by $r_3^0$ (see \eqref{eq:2_order_result_2}),
\begin{equation*}
\begin{aligned}
    \forall \ \left(v_1^*,v_2^*\right)\in \left(V_\omega\right)^2, \quad & \sum_{\alpha,\beta,\gamma,\delta} \int_{\omega} v_{\alpha,\beta}^* \left(\int_\mathcal{I} r_{\alpha}^0\, A_{\alpha\beta\gamma\delta}\, r_{\gamma}^0 \right) v_{\gamma,\delta}^0  \:+\: \sum_{\alpha,\beta} \int_\omega v_{\alpha,\beta}^* \left(\int_\mathcal{I} r_\alpha^0\, A_{\alpha\beta33}\,r_{3,3}^2\right)v_3^0  \\
    + & \sum_{\alpha,\gamma} \int_{\omega} v_\alpha^*\left[ \left(\int_\mathcal{I} r_{3}^0\, A_{\alpha 3\gamma 3}\, r_{3}^0\right) v_\gamma^2  \:+\: \left(\int_\mathcal{I} r_{3}^0\, A_{\alpha3\gamma3}\,r_{3}^0\right) v_{3,\gamma}^2\right]  \\
    + & \sum_{\alpha,\gamma} \int_{\omega} v_\alpha^*\left[ \left(\int_\mathcal{I} r_{3}^0\, A_{\alpha 3\gamma 3}\, r_{\gamma,3}^2\right) v_\gamma^0  \:+\: \left(\int_\mathcal{I} r_{3}^0\, A_{\alpha3\gamma3}\,r_{3}^2\right) v_{3,\gamma}^0\right]  \\
    + & \sum_{\alpha,\gamma} \int_{\omega} v_\alpha^*\left[ \left(\int_\mathcal{I} r_{3}^0\, A_{\alpha 3\gamma 3}\, r_{\gamma,3}^1\right) v_\gamma^1  \:+\: \left(\int_\mathcal{I} r_{3}^0\, A_{\alpha3\gamma3}\,r_{3}^1\right) v_{3,\gamma}^1\right] \:=\: 0.
\end{aligned}
\end{equation*}
Taking functions $v_\alpha^*\in V_\omega$ of the form $v_\alpha^*=-v_{3,\alpha}^*$ with $v_3^*\in H^2_0(\omega)$ and using the fact that $v_\gamma^0=-v_{3,\gamma}^0$, we obtain
\begin{equation} \label{eq:4_order_equations_1}
\begin{aligned}
    \forall \ v_3^*\in H^2_0(\omega), \quad & \sum_{\alpha,\beta,\gamma,\delta} \int_{\omega} v_{3,\alpha\beta}^* \left(\int_\mathcal{I} r_{\alpha}^0\, A_{\alpha\beta\gamma\delta}\, r_{\gamma}^0 \right) v_{3,\gamma\delta}^0  \:-\: \sum_{\alpha,\beta} \int_\omega v_{3,\alpha\beta}^* \left(\int_\mathcal{I} r_\alpha^0\, A_{\alpha\beta33}\,r_{3,3}^2\right)v_{3}^0 \\
    - & \sum_{\alpha,\gamma} \int_{\omega} v_{3,\alpha}^*\left[ \left(\int_\mathcal{I} r_{3}^0\, A_{\alpha 3\gamma 3}\, r_{3}^0\right) v_\gamma^2  \:+\: \left(\int_\mathcal{I} r_{3}^0\, A_{\alpha3\gamma3}\,r_{3}^0\right) v_{3,\gamma}^2\right]  \\
    - & \sum_{\alpha,\gamma} \int_{\omega} v_{3,\alpha}^*\left[ \left(\int_\mathcal{I} r_{3}^0\, A_{\alpha 3\gamma 3}\, r_{\gamma,3}^2\right) v_\gamma^0  \:+\: \left(\int_\mathcal{I} r_{3}^0\, A_{\alpha3\gamma3}\,r_{3}^2\right) v_{3,\gamma}^0\right] \\
    - & \sum_{\alpha,\gamma} \int_{\omega} v_{3,\alpha}^*\left[ \left(\int_\mathcal{I} r_{3}^0\, A_{\alpha 3\gamma 3}\, r_{\gamma,3}^1\right) v_\gamma^1  \:+\: \left(\int_\mathcal{I} r_{3}^0\, A_{\alpha3\gamma3}\,r_{3}^1\right) v_{3,\gamma}^1\right] \:=\: 0.
\end{aligned}
\end{equation}
Similarly, considering \eqref{eq:PGD_scaled_problem_1} now with a test function $\boldsymbol{v}^*$ such that only its third component does not vanish, and taking again into account \eqref{eq:0_order_result}, \eqref{eq:1_order_result}, \eqref{eq:2_order_result_1}, \eqref{eq:2_order_result_2} and \eqref{eq:3_order_result}, we have
\begin{equation} \label{eq:4_order_equations_2}
\begin{aligned}
    \forall \ v_3^*\in V_\omega, \quad & \sum_{\alpha,\gamma} \int_{\omega} v_{3,\alpha}^* \left[\left(\int_\mathcal{I} r_{3}^0\, A_{\alpha 3\gamma 3}\, r_{3}^0\right) v_\gamma^2 \:+\: \left(\int_\mathcal{I} r_3^0\, A_{\alpha3\gamma3}\,r_3^0\right)v_{3,\gamma}^2 \right]  \\
    + & \sum_{\alpha,\gamma} \int_{\omega} v_{3,\alpha}^* \left[\left(\int_\mathcal{I} r_{3}^0\, A_{\alpha 3\gamma 3}\, r_{\gamma,3}^2\right) v_\gamma^0 \:+\: \left(\int_\mathcal{I} r_3^0\, A_{\alpha3\gamma3}\,r_3^2\right)v_{3,\gamma}^0 \right]   \\
    + & \sum_{\alpha,\gamma} \int_{\omega} v_{3,\alpha}^* \left[\left(\int_\mathcal{I} r_{3}^0\, A_{\alpha 3\gamma 3}\, r_{\gamma,3}^1\right) v_{\gamma}^1 \:+\: \left(\int_\mathcal{I} r_3^0\, A_{\alpha3\gamma3}\,r_3^1\right)v_{3,\gamma}^1 \right]  \\
    + & \sum_{\gamma,\delta} \int_\omega v_3^* \left(\int_\mathcal{I} r_{3,3}^2\, A_{33\gamma\delta}\,r_\gamma^0\right)v_{\gamma,\delta}^0 \:+\:  \int_\omega v_3^* \left(\int_\mathcal{I} r_{3,3}^2\, A_{3333}\,r_{3,3}^2\right)v_{3}^0 \:=\: \int_{\omega} v_3^* \left[\left(\int_\mathcal{I} r_3^0\,f_3\right) \:+\: 2\,r_3^0\,g_3\right].
\end{aligned}
\end{equation}
Restricting \eqref{eq:4_order_equations_2} to $v_3^*\in H^2_0(\omega)$ and summing with \eqref{eq:4_order_equations_1}, we obtain
\begin{equation} \label{eq:4_order_equations_3}
\begin{aligned}
    \forall \ v_3^*\in H^2_0(\omega), \quad & \sum_{\alpha,\beta,\gamma,\delta} \int_{\omega} v_{3,\alpha\beta}^* \left(\int_\mathcal{I} r_{\alpha}^0\, A_{\alpha\beta\gamma\delta}\, r_{\gamma}^0 \right) v_{3,\gamma\delta}^0  \:-\: \sum_{\alpha,\beta} \int_\omega v_{3,\alpha\beta}^* \left(\int_\mathcal{I} r_\alpha^0\, A_{\alpha\beta33}\,r_{3,3}^2\right)v_{3}^0 \\
    + & \sum_{\gamma,\delta} \int_\omega v_3^* \left(\int_\mathcal{I} r_{3,3}^2\, A_{33\gamma\delta}\,r_\gamma^0\right)v_{\gamma,\delta}^0 \:+\:  \int_\omega v_3^* \left(\int_\mathcal{I} r_{3,3}^2\, A_{3333}\,r_{3,3}^2\right)v_{3}^0 \: =\: \int_{\omega} v_3^* \left[\left(\int_\mathcal{I} r_3^0\,f_3\right) \:+\: 2\,r_3^0\,g_3\right]
\end{aligned}
\end{equation}
which can be written in the sense of distributions for functions $v_3^*$ in the space $\mathcal{D}(\omega)$. Using the derivation rules on distributions, we obtain
\begin{equation*} \label{eq:4_order_equations_4}
\begin{aligned}
    \forall \ v_3^*\in \mathcal{D}(\omega), \quad & \sum_{\alpha,\beta,\gamma,\delta} \Biggl \langle  \left(\int_\mathcal{I} r_{\alpha}^0\, A_{\alpha\beta\gamma\delta}\, r_{\gamma}^0 \right) v_{3,\gamma\delta\beta\alpha}^0, v_{3}^* \Biggr \rangle  \:-\: \sum_{\alpha,\beta} \Biggl \langle \left(\int_\mathcal{I} r_\alpha^0\, A_{\alpha\beta33}\,r_{3,3}^2\right)v_{3,\beta\alpha}^0, v_{3}^* \Biggr \rangle \\ 
    + & \sum_{\gamma,\delta} \Biggl \langle  \left(\int_\mathcal{I} r_{3,3}^2\, A_{33\gamma\delta}\,r_\gamma^0\right)v_{\gamma,\delta}^0, v_3^* \Biggr \rangle \:+\:  \Biggl \langle  \left(\int_\mathcal{I} r_{3,3}^2\, A_{3333}\,r_{3,3}^2\right)v_{3}^0, v_3^* \Biggr \rangle \: =\: \Biggl \langle \left(\int_\mathcal{I} r_3^0\,f_3\right) \:+\: 2\,r_3^0\,g_3, v_3^* \Biggr \rangle,
\end{aligned}
\end{equation*}
from which we deduce, in the sense of distributions, that
\begin{multline} \label{eq:4_order_int_3}
    \sum_{\alpha,\beta,\gamma,\delta} \left(\int_\mathcal{I} r_{\alpha}^0\:A_{\alpha\beta\gamma\delta}\: r_{\gamma}^0 \right) v_{3,\gamma\delta\beta\alpha}^0  \:-\: \sum_{\alpha,\beta} \left(\int_\mathcal{I} r_\alpha^0\: A_{\alpha\beta33}\: r_{3,3}^2\right)v_{3,\beta\alpha}^0 \:+\: \sum_{\gamma,\delta}  \left(\int_\mathcal{I} r_{3,3}^2\: A_{33\gamma\delta}\:r_\gamma^0\right)v_{\gamma,\delta}^0  \\ 
    + \left(\int_\mathcal{I} r_{3,3}^2\: A_{3333}\: r_{3,3}^2\right)v_3^0  
    = \left(\int_{\mathcal{I}} r_3^0\: f_3\right) \:+\: 2\,r_3^0\,g_3. 
\end{multline}
Using \eqref{eq:2_order_result_1} and the symmetry of $\boldsymbol{C}$, we compute that
\begin{equation*}
    - \sum_{\alpha,\beta} \left(\int_\mathcal{I} r_\alpha^0\: A_{\alpha\beta33}\: r_{3,3}^2\right)v_{3,\beta\alpha}^0 \:+\: \sum_{\gamma,\delta}  \left(\int_\mathcal{I} r_{3,3}^2\: A_{33\gamma\delta}\:r_\gamma^0\right)v_{\gamma,\delta}^0 = -2\sum_{\alpha,\beta} \left(\int_\mathcal{I} r_\alpha^0\: A_{\alpha\beta33}\: r_{3,3}^2\right)v_{3,\alpha\beta}^0.
\end{equation*}
Inserting this identity in \eqref{eq:4_order_int_3}, taking into account \eqref{eq:2_order_result_2} and \eqref{eq:2_order_result_3}, and denoting $u_3^0 = r_3^0\, v_3^0$, Equation~\eqref{eq:4_order_int_3} becomes
\begin{equation} \label{eq:4_order_result}
\boxed{
    \sum_{\alpha,\beta,\gamma,\delta} \left(\int_\mathcal{I} x_3^2\, A_{\alpha\beta\gamma\delta} \right) u_{3,\gamma\delta\beta\alpha}^0 - \left(\int_\mathcal{I} x_3^2\frac{A_{\alpha\beta33}\,A_{33\gamma\delta}}{A_{3333}}\right)\left[2\frac{\displaystyle \int_\omega u_3^0\, u_{3,\gamma\delta}^0}{\displaystyle \int_\omega \left(u_3^0\right)^2} u_{3,\alpha\beta}^0 - \frac{\displaystyle \int_\omega u_3^0\, u_{3,\alpha\beta}^0\int_\omega u_3^0\, u_{3,\gamma\delta}^0}{\displaystyle \left(\int_\omega \left(u_3^0\right)^2\right)^2} u_3^0 \right] 
    = p_3}
\end{equation}
where $p_3$ is given by $\displaystyle p_3(x_1,x_2) = \int_\mathcal{I} f_3(\boldsymbol{x})\,dx_3 \:+\: 2\,g_3(x_1,x_2)$.
\medskip

According to \eqref{eq:bc_v}, the boundary conditions satisfied by $u_3^0$ are
\begin{equation} \label{eq:bc_clamped}
    \left\{
        \begin{array}{rl}
            u_3^0 &= 0 \\
            \boldsymbol{\nabla} u_3^0 \cdot \boldsymbol{n} &= 0
        \end{array}
    \right. \quad \mathrm{on} \quad  \partial\omega.
\end{equation}

\medskip 

Note that the validity of \eqref{eq:4_order_result} is not restricted to the case of a clamped plate. In fact, we use the absence of boundary conditions on the out-of-plane component $\boldsymbol{r}$ (below \eqref{eq:2_order_int_7}), but the boundary conditions on the in-plane component $\boldsymbol{v}$ (see the decomposition \eqref{eq:PGD_decomp}) did not play a role in the derivation of \eqref{eq:4_order_result} (in constrast to the derivation of \eqref{eq:bc_clamped}). It will be useful in the following to consider a simply supported plate. In the case of soft simple support, the displacement field $\boldsymbol{u}$ belongs to the space 
\begin{equation*}
    V = \left\{\boldsymbol{u}\in \left[H^1(\Omega)\right]^3,\; u_3=0\: \text{on}\: \Gamma_0\right\},
\end{equation*}
while in the case of hard simple support $\boldsymbol{u}$ belongs to 
\begin{equation*}
    V = \left\{\boldsymbol{u}\in \left[H^1(\Omega)\right]^3,\; u_3=0\: \text{on}\: \Gamma_0,\; \boldsymbol{u}\cdot\boldsymbol{t}=0\: \text{on}\: \Gamma_0\right\},
\end{equation*}
where $\boldsymbol{t}$ is the tangent vector to $\partial\omega$. Consequently, $\boldsymbol{r}$ is an element of $\displaystyle V_3=H^1\left(-\frac{1}{2},\frac{1}{2}\right)$ and $\boldsymbol{v}$ belongs to either $\displaystyle \left\{\boldsymbol{v}\in \left(V_\omega\right)^3,\; v_3=0\: \text{on}\: \Gamma_0\right\}$ with $V_\omega=H^1(\omega)$ for a soft simple support or $\displaystyle \left\{\boldsymbol{v}\in \left(V_\omega\right)^3,\; v_3=\boldsymbol{v}\cdot\boldsymbol{t}=0\: \text{on}\: \Gamma_0\right\}$ for a hard simple support. Noting that for any $v_3^*\in H^2(\omega)\cap H^1_0(\omega)$, $\boldsymbol{\nabla}v_3^*\cdot \boldsymbol{t}=0$ on $\partial\omega$, \eqref{eq:4_order_equations_1} is satisfied for any $v_3^*\in H^2(\omega)\cap H^1_0(\omega)$ regardless of the type of simple support. We can also restrict \eqref{eq:4_order_equations_2} to $v_3^*\in H^2(\omega)\cap H^1_0(\omega)$ and we obtain as before the analogue of \eqref{eq:4_order_equations_3}, which is written as 
\begin{equation} \label{eq:bc_ss_int_1}
    \forall\, v_3^*\in H^2(\omega)\cap H^1_0(\omega), \quad \sum_{\alpha,\beta} \int_\omega v_{3,\alpha\beta}^* \,\tau_{\alpha\beta} + \int_\omega v_3^* \, \tau_0 = \int_\omega v_3^* \left[\left(\int_\mathcal{I} r_3^0\,f_3\right) \:+\: 2\,r_3^0\,g_3\right],
\end{equation}
where $\tau_{\alpha\beta}$ involves $v_{3,\gamma\delta}^0$ and $v_3^0$ and $\tau_0$ involves $v_{\gamma,\delta}^0$ and $v_3^0$. At this stage, we know that $v_\gamma^0 = -v_{3,\gamma}^0$ and since $v_\gamma^0\in H^1(\omega)$ and $v_3^0\in H^1_0(\omega)$, we have $v_3^0\in H^2(\omega)\cap H^1_0(\omega)$. This implies that $\tau_{\alpha\beta}\in L^2(\omega)$ and $\tau_0 \in L^2(\omega)$. By decomposing $\boldsymbol{\nabla} v_{3}^*$ in the coordinate system $(\boldsymbol{n},\boldsymbol{t})$ and since the tangential derivative of $v_3^*\in H^2(\omega)\cap H^1_0(\omega)$ vanishes, we can define $\displaystyle \sum_{\alpha\beta} n_\alpha\, n_\beta \, \tau_{\alpha\beta}$ on $\partial\omega$ by duality as follows:
\begin{equation} \label{eq:bc_ss_int_2}
    \forall\, v_3^* \in H^2(\omega)\cap H^1_0(\omega), \quad \Biggl \langle \sum_{\alpha,\beta} n_\alpha\, n_\beta \, \tau_{\alpha\beta}, v_{3,n}^* \Biggr \rangle = \sum_{\alpha,\beta} \int_\omega v_{3,\alpha\beta}^* \, \tau_{\alpha\beta} - \int_\omega v_3^* \sum_{\alpha,\beta} \tau_{\alpha\beta,\alpha\beta}.
\end{equation}
Using \eqref{eq:bc_ss_int_1}, \eqref{eq:bc_ss_int_2} becomes 
\begin{equation*} 
    \forall\, v_3^* \in H^2(\omega)\cap H^1_0(\omega), \quad \Biggl \langle \sum_{\alpha,\beta} n_\alpha\, n_\beta \, \tau_{\alpha\beta}, v_{3,n}^* \Biggr \rangle = \int_\omega v_3^* \left[\left(\int_\mathcal{I} r_3^0\,f_3\right) \:+\: 2\,r_3^0\,g_3\right] - \int_\omega v_3^* \, \tau_0 - \int_\omega v_3^* \sum_{\alpha,\beta} \tau_{\alpha\beta,\alpha\beta}.
\end{equation*}
Writing the variational formulation \eqref{eq:bc_ss_int_1} for functions $v_3^*$ in the space $\mathcal{D}(\omega)$, we deduce, in the sense of distribution, that 
\begin{equation*}
    \sum_{\alpha,\beta} \tau_{\alpha\beta,\alpha\beta} + \tau_0 = \left(\int_\mathcal{I} r_3^0\,f_3\right) \:+\: 2\,r_3^0\,g_3,
\end{equation*}
which implies that 
\begin{equation} \label{eq:bc_ss_int_3}
    \forall\, v_3^* \in H^2(\omega)\cap H^1_0(\omega), \quad \Biggl \langle \sum_{\alpha,\beta} n_\alpha\, n_\beta \, \tau_{\alpha\beta}, v_{3,n}^* \Biggr \rangle = 0.
\end{equation}
Regardless of the type of simple support, $v_{3,n}^*$ can take arbitrary values on $\partial\omega$. We thus deduce from \eqref{eq:bc_ss_int_3}, by replacing $\tau_{\alpha\beta}$ with its expression and recalling that $v_3^0=0$ on $\partial\omega$, that
\begin{equation*}
    \sum_{\alpha,\beta} n_\alpha\, n_\beta \sum_{\gamma,\delta} \left(\int_\mathcal{I} r_\alpha^0\, A_{\alpha\beta\gamma\delta}\, r_\gamma^0 \right) v_{3,\gamma\delta}^0  = 0 \quad \mathrm{on} \quad \partial\omega, 
\end{equation*}
that becomes, using \eqref{eq:2_order_result_2}, 
\begin{equation*}
    \sum_{\alpha,\beta} n_\alpha\, n_\beta \sum_{\gamma,\delta}\left(\int_\mathcal{I} x_3^2\, A_{\alpha\beta\gamma\delta} \right) u_{3,\gamma\delta}^0 = 0 \quad \mathrm{on} \quad \partial\omega. 
\end{equation*}
In the case of simple support, the boundary conditions associated with \eqref{eq:4_order_result} are therefore 
\begin{equation} \label{eq:bc_ss} 
    \left\{
        \begin{array}{rl}
            u_3^0 &= 0 \\
            \displaystyle \sum_{\alpha,\beta} n_\alpha\, n_\beta \sum_{\gamma,\delta}\left(\int_\mathcal{I} x_3^2\, A_{\alpha\beta\gamma\delta} \right) u_{3,\gamma\delta}^0 &= 0
        \end{array}
    \right. \quad \mathrm{on} \quad  \partial\omega.
\end{equation}

\subsection{Comparison with the Kirchhoff-Love model: the case of dimension 2} \label{sec:analysis_strip}

\begin{figure}[h!]
    \centering
    \includegraphics[trim=2cm 22cm 2cm 2cm, clip=true, width=\textwidth]{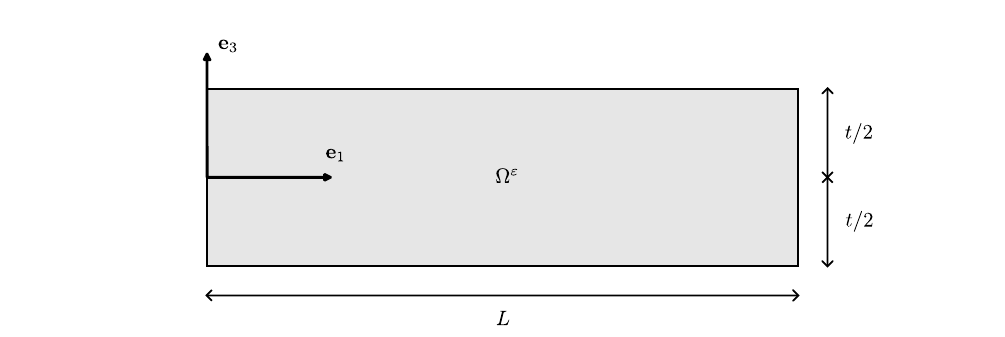}
    \caption{Strip geometry in dimension 2}
    \label{fig:strip_geo}
\end{figure}

To simplify the analysis of the previous results, let us first consider a slightly different problem in which the plate is in a plane strain state with respect to the $(x_1^\varepsilon,x_3^\varepsilon)$ plane. In this subsection, we therefore consider a bending strip of length $L$ and height $t$ in Cartesian coordinates $(x_1^\varepsilon,x_3^\varepsilon)$, where the $x_1^\varepsilon$ axis is along the midline of the strip, while the $x_3^\varepsilon$ axis is along the thickness (see Figure \ref{fig:strip_geo}). The strip therefore occupies the domain $\Omega^{\varepsilon}$ with boundaries $\Gamma_0^{\varepsilon}$ and $\Gamma_{\pm}^{\varepsilon}$ defined by \eqref{eq:plate_domain} where $\omega^L$ is simply the interval $(0,L)$ and $\partial\omega^L$ the set $\{0,L\}$.

This plane strain linear elasticity problem is equivalent to imposing that $u_2=0$ and $u_{i,2}=0$ in all of the above.  
Considering only the first mode of the PGD, the displacement field then takes the form
\begin{equation*} 
    \boldsymbol{u}(x_1,x_3) = \boldsymbol{r}(x_3) \circ \boldsymbol{v}(x_1) = 
    \begin{pmatrix}
        r_1(x_3)\, v_1(x_1) \\
        r_3(x_3)\, v_3(x_1) 
    \end{pmatrix},
\end{equation*}
where $\displaystyle (v_i,r_i)\in H^1(0,1) \times H^1\left(-\frac{1}{2},\frac{1}{2}\right)$ for $i\in\{1,3\}$. The asymptotic development has shown that $r_3^0$ is a constant, that $r_1^0$ is given by
\begin{equation*}
    r_1^0(x_3) = r_3^0\, x_3
\end{equation*}
and 
\begin{equation*}
    v_1^0 = -v_{3,1}^0,
\end{equation*}
so that $v_3^0\in H^2(0,1)$ and, at the leading order,
\begin{equation*} \label{eq:strip_kin}
    \boldsymbol{u}^0(x_1,x_3) = \boldsymbol{r}^0(x_3) \circ \boldsymbol{v}^0(x_1) = 
    \begin{pmatrix}
        -x_3\, r_3^0\, v_{3,1}^0(x_1) \\
        r_3^0\, v_3^0(x_1)
    \end{pmatrix} =
    \begin{pmatrix}
        -x_3\, u_{3,1}^0(x_1) \\
        u_3^0(x_1)
    \end{pmatrix}.
\end{equation*}
As for \eqref{eq:KL_displ}, the kinematics of the first PGD mode thus coincides in the asymptotic limit with a Kirchhoff-Love type kinematics. It remains to establish whether the equilibrium equation and the boundary conditions satisfied by the deflection $w = r_3^0\, v_3^0$ correspond to that of this model. In the 2D case, we find that $w\in H^2(0,1)$ is a solution to
\begin{equation} \label{eq:PGD_strip_equilibrium}
    \left(\int_{\mathcal{I}} x_3^2\, A_{1111} \right) w^{(4)} - \left(\int_{\mathcal{I}}x_3^2 \frac{A_{1133}\,A_{3311}}{A_{3333}}\right) \left[ 2\frac{\displaystyle \int_0^1 w\, w''}{\displaystyle \int_0^1 w^2} w'' - \left(\frac{\displaystyle \int_0^1 w\, w''}{\displaystyle \int_0^1 w^2}\right)^2 w \right]  = p_3
\end{equation}
where the prime denotes the derivative in the distribution sense (with $p_3$ again defined as below \eqref{eq:4_order_result}). For a strip clamped at both ends, the boundary conditions \eqref{eq:bc_clamped} are written as
\begin{equation} \label{eq:bc_strip_cc}
   \left\{
        \begin{array}{rl}
            w\, (0)  = w\,(1) & = 0 \\
            w'(0)  = w'(1) & = 0.
        \end{array}
    \right.
\end{equation}
In the case of a simple support, the boundary conditions \eqref{eq:bc_ss} become
\begin{equation} \label{eq:bc_strip_ss}
    \left\{
        \begin{array}{rl}
            w\,(0) = w\,(1) & = 0 \\
            w''(0) = w''(1) & = 0.
        \end{array}
    \right.
\end{equation}

\subsubsection{The isotropic homogeneous case}

For a plate made of a homogeneous, isotropic material with Young modulus $E$ and Poisson ratio $\nu$, the equilibrium equation \eqref{eq:PGD_strip_equilibrium} becomes
\begin{equation}  \label{eq:PGD_strip_eq_hom}
    \frac{(1-\nu)}{12(1+\nu)(1-2\nu)} w^{(4)} - \frac{\nu^2}{12(1-\nu^2)(1-2\nu)} \left[ 2\frac{\displaystyle \int_0^1 w\, w''}{\displaystyle \int_0^1 w^2} w'' - \left(\frac{\displaystyle \int_0^1 w\, w''}{\displaystyle \int_0^1 w^2}\right)^2 w \right]  = p_3
\end{equation}
where $E$ was chosen as the reference stiffness. Using the above notations, we recall that the scaled equilibrium equation of a bending strip in the Kirchhoff-Love theory is
\begin{equation} \label{eq:bending_strip_equilibrium}
\begin{aligned}
    \frac{1}{12(1-\nu^2)} w^{(4)} = p_3,
\end{aligned} 
\end{equation}
with the same boundary conditions \eqref{eq:bc_strip_cc} or \eqref{eq:bc_strip_ss}. The question that arises is whether a solution to \eqref{eq:bending_strip_equilibrium} is also a solution to \eqref{eq:PGD_strip_eq_hom}. Subtracting Equation~\eqref{eq:bending_strip_equilibrium} from \eqref{eq:PGD_strip_eq_hom}, a necessary condition is that 
\begin{equation*}
    \frac{\nu^2}{12(1-\nu^2)(1-2\nu)}\left[w^{(4)} - 2\frac{\displaystyle \int_0^1 w\, w''}{\displaystyle \int_0^1 w^2} w'' + \left(\frac{\displaystyle \int_0^1 w\, w''}{\displaystyle \int_0^1 w^2}\right)^2 w \right] = 0.
\end{equation*}
This condition is satisfied if $\nu=0$ and in this case Equations~\eqref{eq:PGD_strip_eq_hom} and \eqref{eq:bending_strip_equilibrium} are indeed identical. For $\nu\ne 0$, we necessarily have
\begin{equation} \label{eq:necessary_condition_1}
    w^{(4)} - 2\frac{\displaystyle \int_0^1 w\, w''}{\displaystyle \int_0^1 w^2} w'' + \left(\frac{\displaystyle \int_0^1 w\, w''}{\displaystyle \int_0^1 w^2}\right)^2 w = 0.
\end{equation}
In view of \eqref{eq:necessary_condition_1} and the fact that $w\in H^2(0,1)$, we get $w^{(4)}\in L^2(0,1)$ and thus $w\in H^4(0,1)$. 
Multiplying \eqref{eq:necessary_condition_1} by $w$ and integrating over $(0,1)$, we get
\begin{equation} \label{eq:necessary_condition_2}
    \left(\int_0^1 w\,w''\right)^2 = \int_0^1 w^2 \int_0^1 w\,w^{(4)}.
\end{equation}
By integrating the last integral of \eqref{eq:necessary_condition_2} by parts twice, we obtain
\begin{equation*} 
    \left(\int_0^1 w\,w''\right)^2 = \left(w\,(1\,)\,w^{(3)}(1)-w\,(0)\,w^{(3)}(0) - w'(1)\,w''(1)+w'(0)\,w''(0) + \int_0^1 (w'')^2\right)\int_0^1 w^2.
\end{equation*}
Whether the strip is clamped or simply supported, we deduce from \eqref{eq:bc_strip_cc}, \eqref{eq:bc_strip_ss} and the above equation that
\begin{equation} \label{eq:Cauchy-Schwarz}
    \left(\int_0^1 w\,w''\right)^2 = \int_0^1 w^2 \int_0^1 (w'')^2.
\end{equation}
This corresponds to equality in the Cauchy-Schwarz inequality. We deduce from \eqref{eq:Cauchy-Schwarz} that there exists a constant $\lambda \in\mathbb{R}$ such that $w'' =\lambda\, w$. For the alternative $\lambda \geq 0$, $w=0$ is the only solution with the boundary conditions $w\,(0)=w\,(1)=0$. We therefore necessarily have $\lambda < 0$ and consequently $w\,(x_1) =  w_n \sin(n\pi x_1)$ with $n\in \mathbb{N}^*$ and $w_n\in\mathbb{R}^*$. The form of this solution can only be valid for a simple support, given that $w'(0)\ne 0$, $w'(1)\ne 0$ and $w''(0)=w''(1)=0$. In view of \eqref{eq:bending_strip_equilibrium}, we see that the loading should be of the form $p_3(x_1) = p_n \sin(n\pi x_1)$. Conversely, for a sinusoidal loading of form $p_3(x_1) = p_n \sin(n\pi x_1)$ and a simply supported strip, the solution to \eqref{eq:bc_strip_ss}-\eqref{eq:bending_strip_equilibrium} is $\displaystyle w\,(x_1)=  \frac{12\,(1-\nu^2)p_n }{(n\pi)^4} \sin(n\pi x_1)$, which is also a solution to \eqref{eq:bc_strip_ss}-\eqref{eq:PGD_strip_eq_hom}. 

\begin{remark}
Although $\displaystyle w\,(x_1) =  \sum_{n\in\mathbb{N}^*} w_n \sin(n\pi x_1)$ is a solution to \eqref{eq:bc_strip_ss}-\eqref{eq:bending_strip_equilibrium} for $\displaystyle p_3(x_1)=\sum_{n\in\mathbb{N}^*} p_n \sin(n\pi x_1)$ with $\displaystyle p_n=\frac{n^4\pi^4}{12(1-\nu^2)}w_n$, it is not a solution to \eqref{eq:PGD_strip_eq_hom}. The nonlinearity of \eqref{eq:PGD_strip_eq_hom} is reminiscent of the nonlinearity of the PGD problem \eqref{eq:PGD_scaled_problem}.
\end{remark}

\subsubsection{The laminated case}

It is worth noting that the asymptotic expansion procedure was carried out under the assumption that the tensor $\boldsymbol{C}$ does not depend on $(x_1^\varepsilon,x_2^\varepsilon)$ and is an even function of $x_3^\varepsilon$. The case of laminated plates (or strips) with a symmetrical stacking sequence falls within that setting. It is therefore interesting to establish a link between the solution given by the first mode of PGD and the solution obtained by the Classical Laminated Plates Theory (CPLT). In the isotropic case, the two kinematics coincide as well as the boundary conditions. It thus remains to compare the equilibrium equations.

For the case of bending considered here, the governing equation of CPLT reduces to (see e.g.~\cite{pagano_exact_1969})
\begin{equation*}
    \left(D-\frac{B^2}{A}\right)(w^\varepsilon)^{(4)} = p_3^\varepsilon
\end{equation*}
where $w^\varepsilon$ is the plate deflection. The parameters $A$, $B$, and $D$ are given by
\begin{equation*}
    A = \int_{-\frac{t}{2}}^{\frac{t}{2}} Q_{1111}, \quad B = \int_{-\frac{t}{2}}^{\frac{t}{2}} x_3^\varepsilon \, Q_{1111}, \quad D = \int_{-\frac{t}{2}}^{\frac{t}{2}} \left(x_3^\varepsilon\right)^2 \,Q_{1111},
\end{equation*}
where $Q$ is the plane stress Hooke tensor with entries $Q_{ijkl}$ defined by
\begin{equation*}
    Q_{ijkl} = C_{ijkl} - \frac{C_{ij33}C_{33kl}}{C_{3333}},
\end{equation*}
which only depend (and are even functions of) $x_3^\varepsilon$.
For a symmetrical stacking sequence, the coefficient $B$ vanishes and the equilibrium equation becomes
\begin{equation} \label{eq:CPLT_sym_eq}
    \left(\int_{-\frac{t}{2}}^{\frac{t}{2}} \left(x_3^\varepsilon\right)^2\, Q_{1111} \right)(w^\varepsilon)^{(4)} = p_3^\varepsilon.
\end{equation}
By scaling \eqref{eq:CPLT_sym_eq}, we obtain
\begin{equation} \label{eq:scaled_CPLT_sym_eq}
    \left(\int_{\mathcal{I}} x_3^2\, A_{1111} -\int_{\mathcal{I}}x_3^2 \frac{A_{1133}\,A_{3311}}{A_{3333}}  \right)w^{(4)} = p_3.
\end{equation}
By subtracting \eqref{eq:scaled_CPLT_sym_eq} from \eqref{eq:PGD_strip_equilibrium}, a necessary condition for a solution $w$ to \eqref{eq:scaled_CPLT_sym_eq} to be a solution to \eqref{eq:PGD_strip_equilibrium} is that
\begin{equation*}
    \left(\int_{\mathcal{I}}x_3^2 \frac{A_{1133}\,A_{3311}}{A_{3333}}\right) \left[w^{(4)} - 2\frac{\displaystyle \int_0^1 w\, w''}{\displaystyle \int_0^1 w^2} w'' + \left(\frac{\displaystyle \int_0^1 w\, w''}{\displaystyle \int_0^1 w^2}\right)^2 w \right] = 0,
\end{equation*}
hence, provided that $A_{1133}$ is non-zero,
\begin{equation*}
    w^{(4)} - 2\frac{\displaystyle \int_0^1 w\, w''}{\displaystyle \int_0^1 w^2} w'' + \left(\frac{\displaystyle \int_0^1 w\, w''}{\displaystyle \int_0^1 w^2}\right)^2 w  = 0.
\end{equation*}
This is the same condition as in the isotropic homogeneous case (see \eqref{eq:necessary_condition_1}) and we can conclude in the same way. 

\subsubsection{Asymptotic inconsistency of the first PGD mode} \label{sec:first_conclusions}

For a 2D bending strip, the above analysis shows that the first mode of the PGD method in general does not correspond to the Kirchhoff-Love limit model, even though it exhibits the correct kinematics. It only matches for specific cases of boundary conditions and loadings. 

\subsection{Back to the three-dimensional problem and conclusion}

Let us return to the analysis of the 3D plate problem and start with the isotropic homogeneous case. Denoting $\Delta$ and $\Delta^2$ the Laplacian and bi-Laplacian operators respectively, Equation~\eqref{eq:4_order_result} is then written as
\begin{equation} \label{eq:PGD_plate_equation_hom}
    \frac{(1-\nu)}{12\,(1+\nu)(1-2\nu)}\Delta^2 u_{3}^0 - \frac{\nu^2}{12\,(1-\nu^2)(1-2\nu)}\left[2\frac{\displaystyle \int_\omega u_3^0\, \Delta u_{3}^0}{\displaystyle \int_\omega \left(u_3^0\right)^2} \Delta u_{3}^0 - \left(\frac{\displaystyle \int_\omega u_3^0\, \Delta u_{3}^0}{\displaystyle \int_\omega \left(u_3^0\right)^2} \right)^2 u_3^0 \right] 
    = p_3,
\end{equation}
where $E$ was chosen as the reference stiffness. Considering first a clamped plate, the boundary conditions are given by \eqref{eq:bc_clamped}.
On the other hand, the scaled Kirchhoff-Love equilibrium equation of a bending plate is 
\begin{equation} \label{eq:bending_plate_equation_hom}
    \frac{1}{12\,(1-\nu^2)}\Delta^2u_{3}^0 = p_3
\end{equation}
with the same boundary conditions \eqref{eq:bc_clamped}.
As before, the question is whether a solution to \eqref{eq:PGD_plate_equation_hom} is also a solution to \eqref{eq:bending_plate_equation_hom}, and if so, under which conditions. By subtracting \eqref{eq:bending_plate_equation_hom} from \eqref{eq:PGD_plate_equation_hom}, a necessary condition is that
\begin{equation*}
    \frac{\nu^2}{12\,(1-\nu^2)(1-2\nu)}\left[\Delta^2u_{3}^0 - 2\frac{\displaystyle \int_\omega u_3^0\, \Delta u_{3}^0}{\displaystyle \int_\omega \left(u_3^0\right)^2} \Delta u_{3}^0 + \left(\frac{\displaystyle \int_\omega u_3^0\, \Delta u_{3}^0}{\displaystyle \int_\omega \left(u_3^0\right)^2}\right)^2 u_3^0 \right] = 0.
\end{equation*}
This condition is satisfied if $\nu=0$. In this case, Equations \eqref{eq:PGD_plate_equation_hom} and \eqref{eq:bending_plate_equation_hom} are identical. Otherwise, $u_3^0$ should satisfy 
\begin{equation} \label{eq:necessary_condition_3}
    \Delta^2 u_{3}^0 - 2\frac{\displaystyle \int_\omega u_3^0\, \Delta u_{3}^0}{\displaystyle \int_\omega \left(u_3^0\right)^2} \Delta u_{3}^0 + \left(\frac{\displaystyle \int_\omega u_3^0\, \Delta u_{3}^0}{\displaystyle \int_\omega \left(u_3^0\right)^2}\right)^2 u_3^0 = 0.
\end{equation}
In view of \eqref{eq:necessary_condition_3} and the fact that $u_3^0\in H^2_0(\omega)$, we get $\Delta^2 u_3^0\in L^2(\omega)$. By multiplying \eqref{eq:necessary_condition_3} by $u_3^0$ and integrating over $\omega$, this condition yields 
\begin{equation} \label{eq:necessary_condition_5}
    \left(\int_\omega u_3^0\,\Delta u_3^0\right)^2 = \int_\omega \left(u_3^0\right)^2 \int_\omega u_3^0\, \Delta^2 u_3^0.
\end{equation}
Using the fact that $u_3^0\in H^2_0(\omega)$, $\Delta^2 u \in L^2(\omega)$ and the density of $C^\infty_c(\omega)$ in $H^2_0(\omega)$ for the $H^2$ norm, we obtain that $\displaystyle \int_\omega u_3^0\, \Delta^2 u_3^0 = \int_\omega \left(\Delta u_3^0 \right)^2$. The necessary condition \eqref{eq:necessary_condition_5} thus becomes 
\begin{equation} \label{eq:necessary_condition_6}
    \left(\int_\omega u_3^0\,\Delta u_3^0\right)^2 = \int_\omega \left(u_3^0\right)^2 \int_\omega \left(\Delta u_3^0\right)^2.
\end{equation}
From the Cauchy-Schwarz inequality, \eqref{eq:necessary_condition_6} requires the existence of some $\lambda\in\mathbb{R}$ such that $\Delta u_3^0 = \lambda\, u_3^0$. Together with \eqref{eq:bc_clamped}, this implies that $u_3^0$ vanishes in $\omega$, which is not consistent with \eqref{eq:PGD_plate_equation_hom}, as soon as $p_3$ is not identically zero. Given the boundary conditions, a solution to \eqref{eq:PGD_plate_equation_hom} cannot therefore be a solution to \eqref{eq:bending_plate_equation_hom}.

Consider now a simply supported plate. The boundary conditions associated with \eqref{eq:PGD_plate_equation_hom} are given by \eqref{eq:bc_ss}, that here becomes
\begin{equation*}
    \left\{
        \begin{array}{rl}
            u_3^0 &= 0 \\
            \displaystyle \sum_{\alpha,\beta}\frac{1}{12\,(1-\nu^2)}\left(\frac{\nu\, (1-\nu)}{(1-2\nu)} \Delta u_3^0\, \delta_{\alpha\beta} + (1-\nu)\,u_{3,\alpha\beta}^0\right) n_\beta\, n_\alpha &= 0
        \end{array}
    \right. \quad \mathrm{on} \quad  \partial\omega,
\end{equation*}
while the boundary conditions of the Kirchhoff-Love model are 
\begin{equation*}
    \left\{
        \begin{array}{rl}
            u_3^0 &= 0 \\
            \displaystyle \sum_{\alpha,\beta}\frac{1}{12\,(1-\nu^2)}\left(\nu\, \Delta u_3^0\, \delta_{\alpha\beta} + (1-\nu)\,u_{3,\alpha\beta}^0\right) n_\beta\, n_\alpha &= 0
        \end{array}
    \right. \quad \mathrm{on} \quad  \partial\omega.
\end{equation*}
If $\nu=0$, these boundary conditions are identical, as the two equilibrium equations. If $\nu\ne 0$, we can only expect the two solutions to coincide if the boundary conditions are the same. This implies that $\Delta u_3^0 = 0$ on $\partial\omega$. Going back to \eqref{eq:necessary_condition_5}, we then also find \eqref{eq:necessary_condition_6} by integrating by parts which implies that $\Delta u_3^0=\lambda\, u_3^0$ for some $\lambda$. Conversely, if there exists $\lambda\in\mathbb{R}$ such that $\Delta u_3^0=\lambda\, u_3^0$, we have $\Delta u_3^0=0$ on $\partial\omega$ (since $u_3^0=0$ on $\partial\omega$) and \eqref{eq:PGD_plate_equation_hom} and \eqref{eq:bending_plate_equation_hom} are identical. As an example, for a rectangular plate subjected to a sinusoidal loading of the form $p_3(x_1,x_2)=p_n\sin(n\pi\,x_1)\sin(n\pi\,x_2)$ with $n\in\mathbb{N}^*$, the solution $\displaystyle u_3^0(x_1,x_2) = \frac{3p_n(1-\nu^2)}{n^4\pi^4}\sin(n\pi\,x_1)\sin(n\pi\,x_2)$ to \eqref{eq:bending_plate_equation_hom} indeed satisfies a Helmholtz equation and is a solution to \eqref{eq:PGD_plate_equation_hom}. However, this solution is only valid for simple support boundary conditions.

\medskip 

We now do not restrict ourserlves to the isotropic case. Under the most general assumptions made here and in particular the fact that the components of $\boldsymbol{A}$ are even functions of $x_3$, the scaled governing equation of a bending plate in the Kirchhoff-Love theory is
\begin{equation} \label{eq:bending_plate_equation}
    \sum_{\alpha,\beta,\gamma,\delta} \left(\int_{\mathcal{I}} x_3^2\, A_{\alpha\beta\gamma\delta} -\int_{\mathcal{I}}x_3^2 \frac{A_{\alpha\beta33}\,A_{33\gamma\delta}}{A_{3333}}  \right)u_{3,\alpha\beta\gamma\delta}^0 = p_3.
\end{equation}
Following the same developments as above, a necessary condition for a solution to \eqref{eq:4_order_result} to be a solution to \eqref{eq:bending_plate_equation} is that
\begin{equation} \label{eq:necessary_condition_gen}
    \sum_{\alpha,\beta,\gamma,\delta} \left(\int_{\mathcal{I}}x_3^2 \frac{A_{\alpha\beta33}\,A_{33\gamma\delta}}{A_{3333}}  \right) \left[u_{3,\alpha\beta\gamma\delta}^0 - 2\frac{\displaystyle \int_\omega u_3^0\, u_{3,\gamma\delta}^0}{\displaystyle \int_\omega \left(u_3^0\right)^2} u_{3,\alpha\beta}^0 + \frac{\displaystyle \int_\omega u_3^0\, u_{3,\alpha\beta}^0\int_\omega u_3^0\, u_{3,\gamma\delta}^0}{\displaystyle \left(\int_\omega \left(u_3^0\right)^2\right)^2} u_3^0 \right] = 0.
\end{equation}
For \eqref{eq:necessary_condition_gen} to be satisfied, it suffices that there exists $\lambda\in \mathbb{R}$ such that, for any $\alpha$ and $\beta$, $\displaystyle \sum_{\gamma,\delta} \left(\int_{\mathcal{I}}x_3^2 \frac{A_{\alpha\beta33}\,A_{33\gamma\delta}}{A_{3333}} \right)u_{3,\gamma\delta}^0= \lambda\, u_3^0$. As for the isotropic homogeneous case, this implies that $u_3^0$ vanishes if we consider a clamped plate. With regard to a simply-supported plate, the boundary conditions the Kirchhoff-Love model are
\begin{equation*}
    \left\{
        \begin{array}{rl}
            u_3^0 &= 0 \\
            \displaystyle \sum_{\alpha,\beta} n_\alpha\, n_\beta \sum_{\gamma,\delta} \left(\int_I x_3^2\left(A_{\alpha\beta\gamma\delta}-\frac{A_{\alpha\beta33}A_{33\gamma\delta}}{A_{3333}}\right) \right)u_{3,\gamma\delta}^0  &= 0
        \end{array}
    \right. \quad \mathrm{on} \quad  \partial\omega.
\end{equation*}
The above conditions and \eqref{eq:bc_ss} can only be identical if, for any $\alpha$ and $\beta$, $\displaystyle \sum_{\gamma,\delta} \left(\int_{\mathcal{I}}x_3^2 \frac{A_{\alpha\beta33}\,A_{33\gamma\delta}}{A_{3333}}  \right)u_{3,\gamma\delta}^0= 0$ on $\partial\omega$, which is consistent with the previous sufficient condition.
\medskip

The method of formal asymptotic expansion has thus enabled us to study the behavior of the displacement field given by the first PGD mode when the plate thickness goes to $0$. Although this method does not constitute a rigorous mathematical proof, the results tend to show that the limit solution given by the first PGD mode is not the asymptotic solution, i.e. the Kirchhoff-Love solution, except in very specific cases of boundary conditions and loadings. The kinematics is the correct one, but the equilibrium equation differs.

\section{Towards a new PGD strategy for slender elastic structures}
\label{sec:new_PGD_strategy}

\subsection{Why does the standard approach fail?}

It has been shown so far that the first mode provided by PGD is not asymptotically consistent. This result is discussed here in the light of other considerations made in the literature. For a homogeneous plate subjected to bending loading, it can be shown according to homogenization theory \cite{lesage_preprint} and under the assumption of symmetry \eqref{eq:monoclinic}, that the 3D elasticity solution, when the thickness goes to 0, expands as follows:
\begin{equation} \label{eq:homogenized_form}
    \boldsymbol{u} (x_1,x_2,x_3) \approx 
    \begin{pmatrix}
        -x_3\, u_{3,1}^\star (x_1,x_2) \\
        -x_3\, u_{3,2}^\star (x_1,x_2) \\
        u_3^\star(x_1,x_2)
    \end{pmatrix} + \varepsilon^2
    \begin{pmatrix}
        0 \\
        0 \\
        \displaystyle \frac{1}{2}\left(x_3^2-\frac{1}{12}\right)\sum_{\alpha,\beta} \frac{A_{33\alpha\beta}}{A_{3333}}u_{3,\alpha\beta}^\star(x_1,x_2)
    \end{pmatrix}.
\end{equation}
The approximation \eqref{eq:homogenized_form}, where $u_3^\star$ is the solution to a 2D homogenized problem, is characterised by a separation of variables, which supports the use of PGD in this context. In addition, the approximation \eqref{eq:homogenized_form} holds in energy norm and not only in $L^2(\Omega)$ norm. Observe that, although the second term in the right hand side of \eqref{eq:homogenized_form} is negligible in terms of displacement, it is significant in terms of energy. This is why, in the energy norm, $\boldsymbol{u}$ is not close to a single PGD mode, but to a sum of two PGD modes. 

On a similar note, it is well known that inserting Kirchhoff-Love kinematics into the 3D elasticity equations does not yield the correct equilibrium equation. The equation obtained is typical of a state of plane strain, whereas the correct assumption is that of plane stress. Finally, in the nomenclature of hierarchical models, it is established that the model $(1,1,0)$ is not asymptotically correct. Interested readers are referred to \cite{braess_justification_2011} or \cite{paumier_asymptotic_1997} for more details on this subject. 

We numerically show in Section~\ref{sec:numerical_results} that the first PGD mode is indeed a poor approximation of the reference solution in the sense that, in the standard PGD strategy (where modes are computed one after each other in an iterative way), this first PGD mode does not converge to the exact solution when the thickness goes to 0. Furthermore, it is challenging to \textit{a priori} assess the number of PGD modes required to obtain an accurate solution.

\subsection{Asymptotic consistency of higher-rank PGD solution}

The question that naturally arises is how many PGD modes are needed to recover asymptotic consistency. Based on the analysis of the first PGD mode, we claim that regardless of the rank of the PGD solution, the latter does not converge to the limit solution when the thickness tends to 0. 

Considering the computation of a second PGD mode of the form $\boldsymbol{s}(x_3)\circ\boldsymbol{w}(x_1,x_2)$, we provide arguments that support this statement. Assuming that the first mode $\boldsymbol{r}\circ\boldsymbol{v}$ is computed, we look for $(\boldsymbol{w},\boldsymbol{s})\in (V_\omega\times V_3)^3$ such that, for any $(\boldsymbol{w}^*,\boldsymbol{s}^*)\in (V_\omega\times V_3)^3$,
\begin{equation*}
\int_{\Omega} \boldsymbol{\epsilon}\,(\boldsymbol{s}\circ \boldsymbol{w}^* + \boldsymbol{s}^*\circ \boldsymbol{w}):\boldsymbol{C}:\boldsymbol{\epsilon}\,(\boldsymbol{r}\circ \boldsymbol{v}+\boldsymbol{s}\circ \boldsymbol{w}) = \int_{\Omega} (\boldsymbol{s}\circ \boldsymbol{w}^* + \boldsymbol{s}^*\circ \boldsymbol{w}) \cdot \boldsymbol{f} + \int_{\Gamma_+\cup\, \Gamma_-} (\boldsymbol{s}\circ \boldsymbol{w}^* + \boldsymbol{s}^*\circ \boldsymbol{w}) \cdot \boldsymbol{g},
\end{equation*}
or even
\begin{equation*}
\int_{\Omega} \boldsymbol{\epsilon}\,(\boldsymbol{s}\circ \boldsymbol{w}^* + \boldsymbol{s}^*\circ \boldsymbol{w}):\boldsymbol{C}:\boldsymbol{\epsilon}\,(\boldsymbol{s}\circ \boldsymbol{w}) = \int_{\Omega} (\boldsymbol{s}\circ \boldsymbol{w}^* + \boldsymbol{s}^*\circ \boldsymbol{w}) \cdot \left(\boldsymbol{f}+\nabla\cdot \boldsymbol{\sigma}\right) + \int_{\Gamma_+\cup\, \Gamma_-} (\boldsymbol{s}\circ \boldsymbol{w}^* + \boldsymbol{s}^*\circ \boldsymbol{w}) \cdot \left(\boldsymbol{g}-\boldsymbol{\sigma}\cdot\boldsymbol{n}\right), 
\end{equation*}
where $\boldsymbol{\sigma}=\boldsymbol{C}:\boldsymbol{\epsilon}\,(\boldsymbol{r}\circ\boldsymbol{v})$. The first PGD mode is thus included in the loading for the computation of the second mode, which reduces to a single-mode problem. In light of the previous analysis, we would expect the second mode not to be asymptotically consistent with respect to the problem with the modified load. There is therefore no \textit{a priori} reason for the sum of the two modes to be asymptotically consistent with respect to the initial problem. This reasoning can be extended to an arbitrary number of modes. 

\begin{remark} \label{rem:asymptotic_consistency}
    We only claim that, for a fixed number of modes $m$, the PGD solution does not converge to the homogenized solution when the thickness tends to 0. In particular, we say nothing about the error between the PGD solution and the exact solution. For a given slenderness ratio that is sufficiently large, the PGD solution converges to the exact solution, which is close to the limit solution. The error between the PGD solution and the limit solution may be small. 
\end{remark}

For all these reasons, we thus present now an alternative PGD strategy, that will be shown to be accurate.

\subsection{Block PGD mode computation}

In view of the above observations, a modification of the standard PGD procedure is proposed, with the aim of capturing the asymptotic solution at the early stage of the procedure. This new PGD strategy, described below, involves computing the first two modes simultaneously, as suggested by the form \eqref{eq:homogenized_form} of the solution provided by homogenization theory. 

\medskip

We consider here the problem in its original formulation \eqref{eq:epsilon_weak_form}, i.e. before scaling. To simplify notation, the exponent $\varepsilon$ is omitted, without risk of confusion. At each iteration $m\in\mathbb{N}^*$, an approximation $\boldsymbol{u}_m\in V$ of the displacement field solution to \eqref{eq:epsilon_weak_form} is constructed in the form
\begin{equation*}
    \boldsymbol{u}_m = \sum_{k=1}^m \boldsymbol{z}_k,
\end{equation*}
where each term of the sum is computed iteratively using a greedy algorithm. The novelty introduced here lies in the modification of the ansatz for the transverse displacement of the first approximation $\boldsymbol{u}_1$. More precisely, $\boldsymbol{u}_1$ is sought in the form
\begin{equation} \label{eq:PGD2_displacement}
    \boldsymbol{u}_1(\boldsymbol{x}) = \boldsymbol{r}(x_3)\circ \boldsymbol{v}(x_1,x_2) + \boldsymbol{s}(x_3)\circ \boldsymbol{w}(x_1,x_2) = 
    \begin{pmatrix}
        r_1(x_3)\, v_1(x_1,x_2) \\
        r_2(x_3)\, v_2(x_1,x_2) \\
        v_3(x_1,x_2) + s_3(x_3)\, w_3(x_1,x_2) 
    \end{pmatrix},
\end{equation}
where $(r_1,r_2,s_3)\in \left(V_3\right)^3$ and $(\boldsymbol{v},w_3)\in \left(V_\omega\right)^4$.

We now minimize the potential energy \eqref{eq:potential_energy} upon displacements of the form \eqref{eq:PGD2_displacement}. Assuming that a minimizer exists, and denoting it by $\boldsymbol{u}_1$, it satisfies the following Euler-Lagrange equation: for any $(\boldsymbol{v}^*,w_3^*)\in \left(V_\omega\right)^4$ and any $(r_1^*,r_2^*,s_3^*)\in \left(V_3\right)^3$,
\begin{equation*}
\begin{aligned}
\int_{\Omega} \boldsymbol{\epsilon}\,(\boldsymbol{r}\circ \boldsymbol{v}^* + \boldsymbol{s}\circ \boldsymbol{w}^* + \boldsymbol{r}^*\circ \boldsymbol{v} + \boldsymbol{s}^*\circ \boldsymbol{w}):\boldsymbol{C}:\boldsymbol{\epsilon}\,(\boldsymbol{u}_1) = \int_{\Omega} (v_3^* + s_3\,w_3^* + s_3^*\, w_3) \, f_3 + \int_{\Gamma_+\cup\, \Gamma_-} (v_3^* + s_3\,w_3^* + s_3^*\, w_3) \, g_3, 
\end{aligned}
\end{equation*}
which can be written equivalently as a system of coupled equations:
\begin{subequations} \label{eq:PGD2_problem}
\begin{align}[left = \empheqlbrace\,]
    &\forall \ (\boldsymbol{v}^*,w_3^*)\in \left(V_\omega\right)^4,  \int_{\Omega} \boldsymbol{\epsilon}\,(\boldsymbol{r}\circ \boldsymbol{v}^* + \boldsymbol{s}\circ \boldsymbol{w}^*):\boldsymbol{C}:\boldsymbol{\epsilon}\,(\boldsymbol{r}\circ \boldsymbol{v} + \boldsymbol{s}\circ \boldsymbol{w}) = \int_{\Omega} (v_3^* + s_3\,w_3^*) \, f_3 + \int_{\Gamma_+\cup\, \Gamma_-} (v_3^* + s_3\,w_3^*) \, g_3 \label{eq:PGD2_problem_1} \\ 
    &\forall \ (r_1^*,r_2^*,s_3^*)\in \left(V_3\right)^3,  \int_{\Omega} \boldsymbol{\epsilon}\,
    (\boldsymbol{r}^*\circ \boldsymbol{v} + \boldsymbol{s}^*\circ \boldsymbol{w}):\boldsymbol{C}:\boldsymbol{\epsilon}\,(\boldsymbol{r}\circ \boldsymbol{v} + \boldsymbol{s}\circ \boldsymbol{w}) = \int_{\Omega} (s_3^*\circ w_3)\, f_3 + \int_{\Gamma_+\cup\, \Gamma_-} (s_3^*\circ w_3) \, g_3.\label{eq:PGD2_problem_2} \
\end{align}
\end{subequations}
In practice, the system of equations \eqref{eq:PGD2_problem} is solved using a fixed-point algorithm. Initial functions $(r_1^{(0)},r_2^{(0)},s_3^{(0)})$ are chosen. Then, at each step $n\geq 1$, the algorithm computes $(\boldsymbol{v}^{(n)},w_3^{(n)},r_1^{(n)},r_2^{(n)},s_3^{(n)})$ such that
\begin{itemize}
    \item $(\boldsymbol{v}^{(n)},w_3^{(n)})$ satisfy equation \eqref{eq:PGD2_problem_1} for $(\boldsymbol{r},\boldsymbol{s})$ set to $(r_1^{(n-1)},r_2^{(n-1)},s_3^{(n-1)})$;
    \item $(r_1^{(n)},r_2^{(n)},s_3^{(n)})$ satisfy equation \eqref{eq:PGD2_problem_2} for $(\boldsymbol{v},\boldsymbol{w})$ set to $(\boldsymbol{v}^{(n)},w_3^{(n)})$.
\end{itemize}
The fixed-point algorithm stops when
\begin{equation} \label{eq:fp_tol}
    \frac{\|\boldsymbol{r}^{(n)}\circ \boldsymbol{v}^{(n)} + \boldsymbol{s}^{(n)}\circ \boldsymbol{w}^{(n)} - (\boldsymbol{r}^{(n-1)}\circ \boldsymbol{v}^{(n-1)} + \boldsymbol{s}^{(n-1)}\circ \boldsymbol{w}^{(n-1)})\|}{\|\boldsymbol{r}^{(n-1)}\circ \boldsymbol{v}^{(n-1)} + \boldsymbol{s}^{(n-1)}\circ \boldsymbol{w}^{(n-1)}\|} < \eta
\end{equation}
where $\eta$ is a predefined tolerance threshold and $\|\cdot \|$ is the energy norm.

Once $\boldsymbol{u}_1$ is known and if necessary, new modes $\boldsymbol{z}_k$ of form $\boldsymbol{z}_k(\boldsymbol{x}) = \boldsymbol{r}_k(x_3)\circ\boldsymbol{v}_k(x_1,x_2)$ can be added to the solution using the standard PGD procedure. 

\begin{remark}
The displacement $\boldsymbol{u}_1$ could have been sought in the form
\begin{equation*}
    \boldsymbol{u}_1(\boldsymbol{x}) = 
    \begin{pmatrix}
        r_1(x_3)\, v_1(x_1,x_2) \\
        r_2(x_3)\, v_2(x_1,x_2) \\
        r_3(x_3)\, v_3(x_1,x_2)
    \end{pmatrix} +
    \begin{pmatrix}
        s_1(x_3)\, w_1(x_1,x_2)  \\
        s_2(x_3)\, w_2(x_1,x_2)  \\
        s_3(x_3)\, w_3(x_1,x_2) 
    \end{pmatrix}.
\end{equation*}
However, in this form, the two sum terms play a similar role and are a priori interchangeable. This may create issues for the convergence of the fixed point iterations. To differentiate between these two terms, the form \eqref{eq:PGD2_displacement} is preferred, which is inspired by the result \eqref{eq:homogenized_form} in pure bending. \\
Under the symmetry assumptions made here on $\boldsymbol{C}$, we expect the PGD approximation \eqref{eq:PGD2_displacement} to be asymptotically consistent for 2D (homogeneous or laminated) strips and for homogeneous plates. However, with regard to laminated plates, the corrector term $s_3\,w_3$ may not be sufficient to recover asymptotic consistency. 
\end{remark}

\section{Numerical experiments}
\label{sec:numerical_results}

In this section, the theoretical conclusions drawn above are numerically validated, and the performance of the proposed new strategy for building a PGD reduced-order model for slender structures is investigated. Section \ref{sec:locking_test} deals with the issue of locking in a PGD context. The inability of the first PGD mode to capture the asymptotic solution is illustrated in Section \ref{sec:first_mode_fail}. The new PGD strategy, based on computing the first two modes simultaneously, is considered in Section \ref{sec:new_PGD_test}. 

All examples are conducted in the representative case of a strip under pure bending, with no body forces and the same density of surface forces applied to the upper and lower faces. The constitutive behavior is assumed to be homogeneous and isotropic, with material parameters
\begin{equation*}
    E = 1~GPa, \quad \nu = 0.3.
\end{equation*}
Different types of boundary conditions and loadings are considered and detailed in Table \ref{tab:BCs_loading}.
\begin{table}[h]
    \centering
    \begin{tabular}{l|c|c}
    \hline
         &  Boundary conditions & Loading \\
    \hline
    SS-SP & Simply Supported & Sinus\\
    SS-UP & Simply Supported & Uniform \\
    CC-SP & Clamped & Sinus \\
    CC-UP & Clamped & Uniform \\
    \hline
    \end{tabular}
    \caption{Boundary conditions and loadings}
    \label{tab:BCs_loading}
\end{table}
Up to an amplitude factor, the sinus loading corresponds to $\displaystyle g_3(x_1) =\sin\left(\frac{\pi x_1}{L}\right)$ and the uniform case to $\displaystyle g_3=1$. Clamped boundary conditions are written as $\boldsymbol{u}(0,x_3)=\boldsymbol{u}(L,x_3)=0$ , while simple support corresponds to $u_3(0,x_3)=u_3(L,x_3)=0$. According to Section~\ref{sec:analysis_strip}, we expect the first PGD mode to be asymptotically correct only in the SS-SP case. 

Unless otherwise specified, 64 quadratic finite elements are used to discretize the axial problem, and a single fourth-order polynomial expansion is used in the thickness. Following the discussion of Section~\ref{sec:locking_test} below, selective reduced integration is used to prevent locking issues and the fixed-point tolerance $\eta$ in \eqref{eq:fp_tol} is set to $10^{-3}$. For more details on implementation, we refer the reader to \ref{app1}. 

A boundary layer may be present in a small area near the boundary depending on the type of boundary conditions. The accuracy of the results obtained depends in part on taking this effect into account. To give the method the possibility to capture this singularity, the finite element mesh  must be able to describe this boundary layer. Although it is possible to use a sufficiently fine uniform mesh, this choice is not optimal in terms of computational cost. Here, we resort a non-uniform mesh pattern, as suggested in the context of hierarchical models \cite{cho_locking_1997}. This mesh consists of a small element of size $0.1\,t$ near the boundary, followed by an element of size $0.9\,t$ and standard elements of size $h$.

\subsection{Demonstrating locking in a PGD context} \label{sec:locking_test}

The phenomenon of shear locking is a common issue in the numerical resolution of plate problems. However, this issue in the PGD framework is only explicitly addressed in a few works including \cite{vidal_robust_2018}. Asymptotic analysis has shown that $v_3^0$ is solution of a variational formulation set in $H^2$. Everything therefore suggests that standard $H^1$- conforming finite elements will lead to the usual shear locking, since they are not comptatible with the regularity of the asymptotic solution. This issue is therefore addressed first, before moving on to a detailed understanding of the asymptotic behavior of PGD modes.

For the sake of clarity, locking in the PGD context is discussed in light of the asymptotic analysis previously conducted for a single mode. The conclusions remain the same when two modes are computed simultaneously, as motivated above. We have observed that the constraint 
\begin{equation} \label{eq:KL_constraint}
    v_\alpha + v_{3,\alpha} = 0
\end{equation}
is enforced when $\varepsilon$ goes to 0 (see \eqref{eq:2_order_result_1}). However, standard low-order finite elements cannot correctly represent this constraint, which causes shear locking. Note that \eqref{eq:KL_constraint} is nothing less than the Kirchhoff-Love constraint (zero shear strain). The bending and transverse shear strains $\boldsymbol{\epsilon}^b$ and $\boldsymbol{\epsilon}^s$ associated with a PGD mode are indeed defined by
\begin{equation*}
    2\, \epsilon_{\alpha\beta}^b = r_\alpha\, v_{\alpha,\beta}+r_\beta\, v_{\beta,\alpha}, \quad 2\,\epsilon_{\alpha 3}^s = r_{\alpha,3}\,v_\alpha + r_3\, v_{3,\alpha}, \quad \epsilon_{33}^b = r_{3,3}\, v_3.
\end{equation*}
The asymptotic limit \eqref{eq:KL_constraint} together with \eqref{eq:2_order_result_2} imposes zero shear strain. Neglecting the out-of-plane normal stress, these strains are associated with the following strain energies $\mathcal{E}^b$ and $\mathcal{E}^s$:
\begin{equation*}
    \mathcal{E}^b = \frac{1}{2}\int_\Omega \epsilon_{\alpha\beta}^b\, C_{\alpha\beta\gamma\delta}\,\epsilon_{\gamma\delta}^b \quad \mathrm{and} \quad \mathcal{E}^s = 2\int_\Omega \epsilon_{\alpha3}^s\, C_{\alpha3\gamma3}\,\epsilon_{\gamma3}^s.
\end{equation*}
According to the variational formulations \eqref{eq:scaled_weak_form} or \eqref{eq:PGD_scaled_problem}, the bending strain energy scales with $\varepsilon^4$ while the transverse shear strain energy scales with $\epsilon^2$. Consequently, if the Kirchhoff-Love constraint cannot be properly represented in the discrete spaces, shear will dominate the total energy and can lead to a significant underestimation of bending strains. 

The purpose of this work is not to review the large number of locking-free numerical methods dedicated to plate models. With regard to the numerical results presented here, a reduced and selective integration method \cite{zienkiewicz_reduced_1971,hughes_simple_1977} is used. The resulting changes in the implementation are presented in \ref{app1}.

\begin{figure}[h!]
    \centering
    \begin{subfigure}{0.49\textwidth}
        \includegraphics[trim=3cm 9cm 4cm 9cm, clip=true, width=\textwidth]{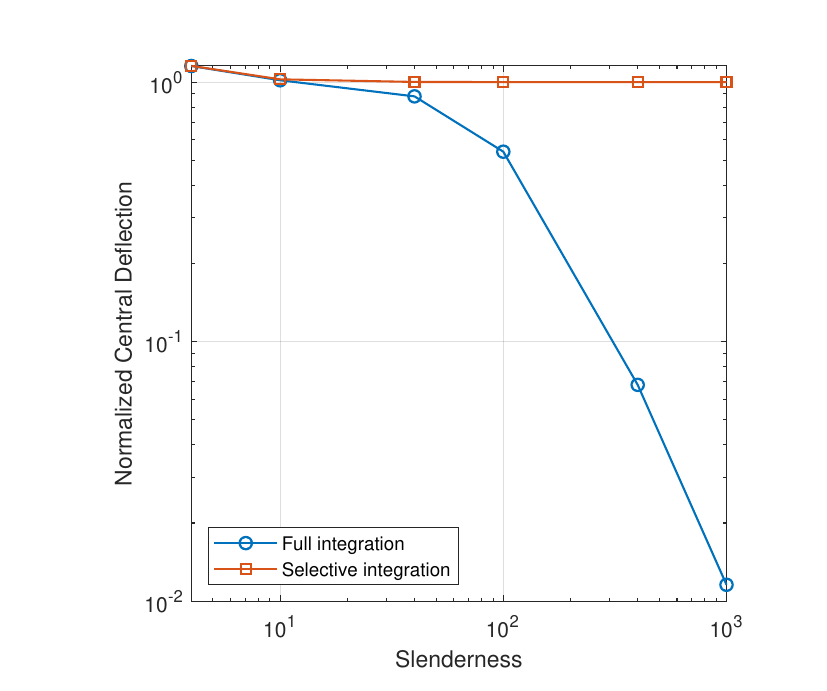}
        \caption{Linear finite elements}
        \label{fig:locking_test_1}
    \end{subfigure}
    \hfill
    \begin{subfigure}{0.49\textwidth}
        \includegraphics[trim=3cm 9cm 4cm 9cm, clip=true, width=\textwidth]{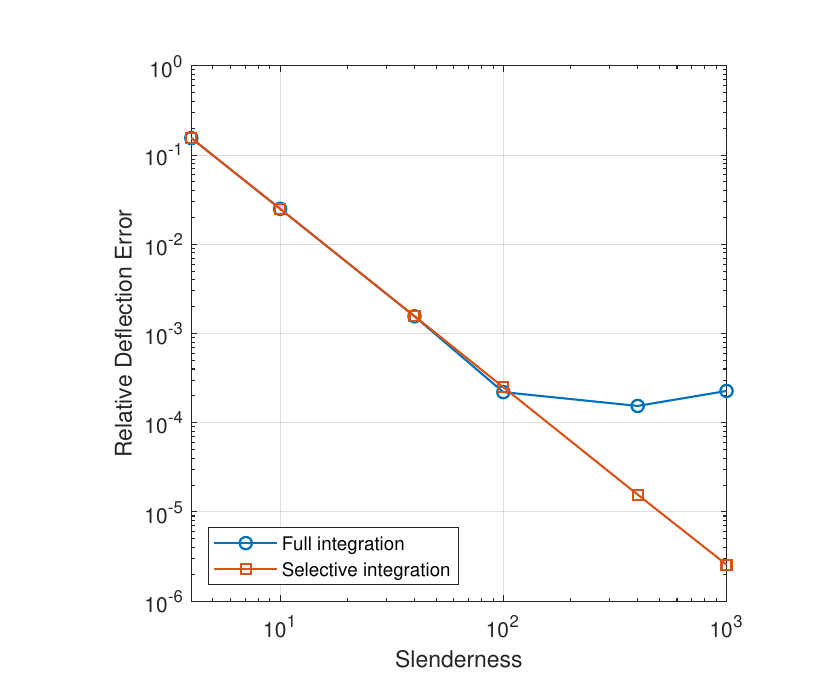}
        \caption{Quadratic finite elements}
        \label{fig:locking_test_2}
    \end{subfigure}
    \caption{Illustration of the transverse shear locking in PGD reduced-order models}
    \label{fig:locking_test}
\end{figure}

To demonstrate that our approach prevents locking, we focus here on the SS-SP case, for which the first PGD mode is asymptotically consistent, i.e. converges (in the absence of any discretization) to the exact solution in the limit of a large slenderness. The deflection at the center of the strip provided by the first PGD mode is shown as a function of the slenderness in Figure \ref{fig:locking_test_1}, where 64 linear finite elements have been used to discretize the axial problem. The deflection is normalized by the deflection value  $w_{KL}$ given by Kirchhoff-Love theory (calculated analytically), and two cases are considered depending on whether selective integration is used or not. The corresponding numerical values are given in Table \ref{tab:locking_test}. It appears that complete shear locking occurs for linear elements (the central deflection goes to zero when slenderness increases, which is completely wrong), and that reduced integration alleviates this problem (the central deflection goes to the asymptotic Kirchhoff-Love value). Figure \ref{fig:locking_test_2} shows the relative deflection error defined by
\begin{equation} \label{eq:relative_defl_err}
    \mathrm{Relative\ Deflection\ Error} = \left\lvert\frac{\displaystyle r_3(0)\,v_3\left(\frac{L}{2}\right) - w_{KL}}{w_{KL}} \right\rvert,
\end{equation}
when quadratic elements are used in the axial direction. Shear locking is considerably less important than in the case of linear elements and appears only for slenderness larger than 100. When using selective integration, we clearly observe convergence of the first PGD mode to the Kirchhoff-Love solution as slenderness increases.

\begin{table}[h!]
    \centering
    \begin{tabular}{lllllllll}
    \hline
    Slenderness: $L/t$ & 4 & 10 & 40 & $10^2$ & $4\times 10^2$ & $10^3$ & $4\times 10^3$ & $10^4$   \\
    \hline 
    Full integration & 1.1529 &  1.0158  &  0.8807   & 0.5391   & 0.0681  &  0.0116 & 0.0007 & 0.0001 \\
    Selective integration & 1.1543  &  1.0243 &   1.0009  &  0.9996 &   0.9994  &  0.9994 & 0.9994 & 0.9994 \\
    \hline
    \end{tabular}
    \caption{Normalized central deflection as a function of slenderness for the SS-SP case and linear finite elements used in the axial direction}
    \label{tab:locking_test}
\end{table}

\subsection{Asymptotic inconsistency of the first PGD mode and higher-rank PGD solutions} \label{sec:first_mode_fail}

The purpose of this section is to show that the first mode provided by PGD is not asymptotically consistent except in certain special cases of loading and boundary conditions. To this end, the convergence of the first PGD mode as a function of slenderness is investigated. Two relative errors are introduced. The first is the relative deflection error \eqref{eq:relative_defl_err} already defined. The second is the relative strain energy error of the first mode in comparison with $\mathcal{E}_{KL}$, the exact strain energy of the displacement field in the Kirchhoff-Love theory (calculated analytically). We recall that $\mathcal{E}_{KL}$ is the limit of the strain energy of the exact solution when the thickness goes to 0. The energy error is defined by
\begin{equation}
    \mathrm{Relative\ Energy\ Error} = \left\lvert\frac{\mathcal{E} - \mathcal{E}_{KL}}{\mathcal{E}_{KL}} \right\rvert,
\end{equation}
where $\mathcal{E}$ is the strain energy of the first PGD mode.

\begin{figure}[h!]
    \centering
    \begin{subfigure}{0.49\textwidth}
        \includegraphics[trim=3cm 9cm 4cm 9cm, clip=true, width=\textwidth]{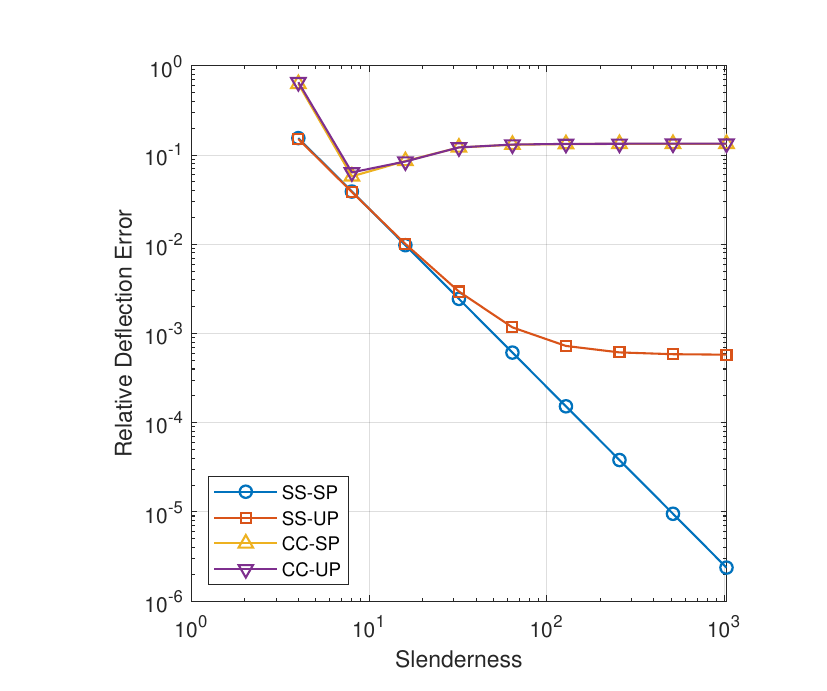}
    \end{subfigure}
    \begin{subfigure}{0.49\textwidth}
        \includegraphics[trim=3cm 9cm 4cm 9cm, clip=true, width=\textwidth]{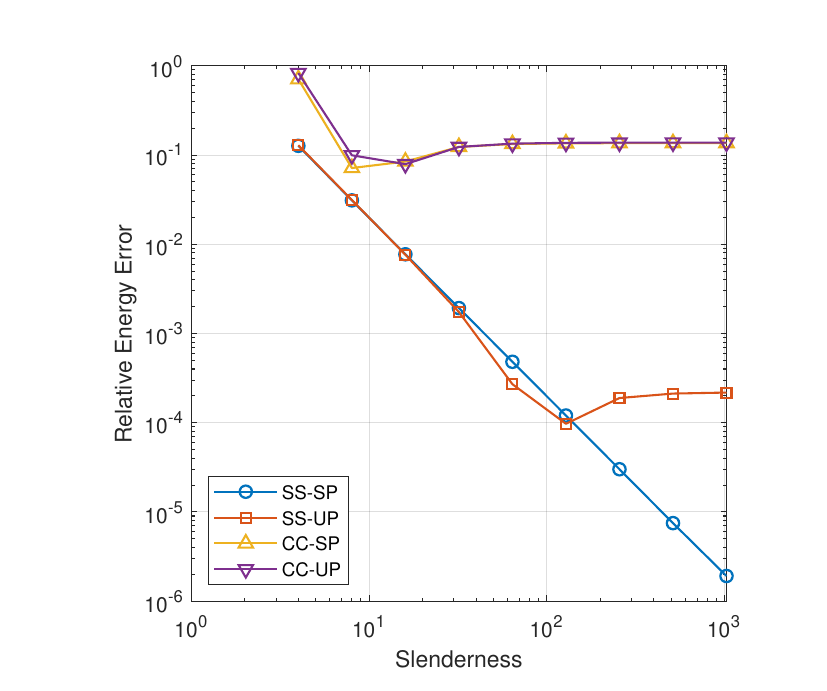}
    \end{subfigure}
    \caption{Relative deflection (left) and energy (right) error of the first PGD mode as a function of slenderness (case $\nu\ne 0$)}
    \label{fig:CV_first_mode}
\end{figure}

\begin{figure}[h!]
    \centering
    \begin{subfigure}{0.49\textwidth}
        \includegraphics[trim=3cm 9cm 4cm 9cm, clip=true, width=\textwidth]{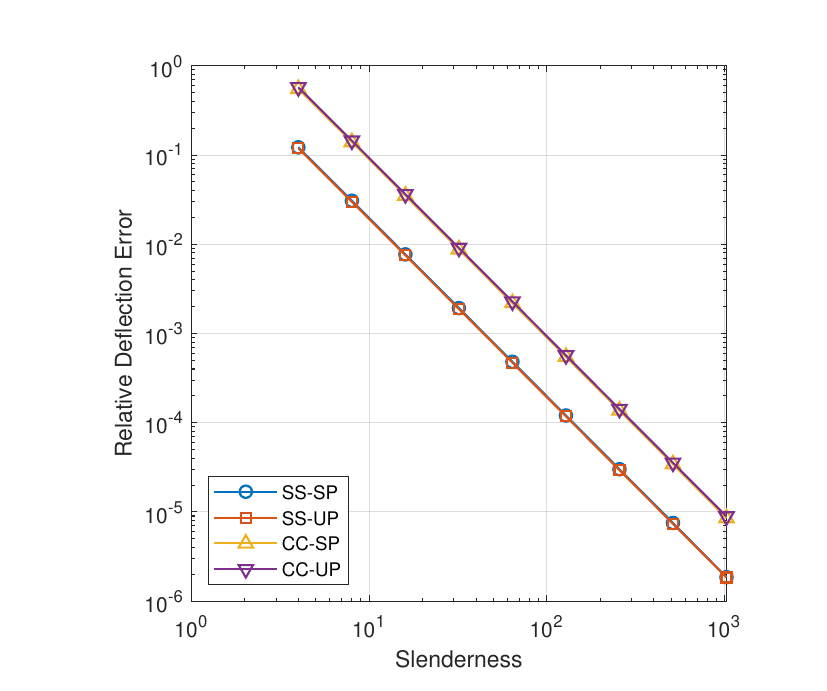}
    \end{subfigure}
    \begin{subfigure}{0.49\textwidth}
        \includegraphics[trim=3cm 9cm 4cm 9cm, clip=true, width=\textwidth]{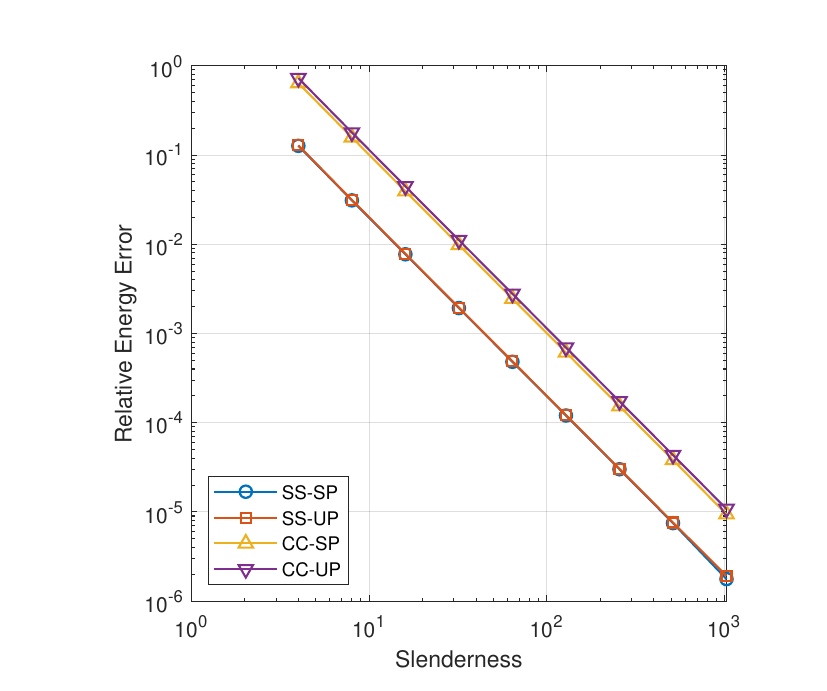}
    \end{subfigure}
    \caption{Relative deflection (left) and energy (right) error of the first PGD mode as a function of slenderness (case $\nu= 0$)}\label{fig:CV_first_mode_nu0}
\end{figure}

These two errors as a function of slenderness are shown in Figure~\ref{fig:CV_first_mode} (for a case where $\nu\ne0$) and in Figure~\ref{fig:CV_first_mode_nu0} (for a case where $\nu=0$, i.e. a case without any Poisson effect). The four cases of boundary conditions and loadings described in Table~\ref{tab:BCs_loading} are presented. Figure~\ref{fig:CV_first_mode} shows that in general the first mode does not converge to the asymptotic solution, except in the SS-SP case corresponding to a simply supported strip subjected to sinusoidal loading. The conclusions drawn in Section \ref{sec:first_conclusions} are thus confirmed numerically. Now considering a zero Poisson ratio (Figure \ref{fig:CV_first_mode_nu0}), the first PGD mode converges to the Kirchhoff-Love solution regardless of the boundary conditions and loadings. This is in agreement with Equation~\eqref{eq:PGD_strip_eq_hom}, which coincides with the equilibrium equation of Kirchhoff-Love theory when $\nu=0$.

\begin{figure}[h!]
    \centering
    \begin{subfigure}{0.49\textwidth}
        \includegraphics[trim=3cm 9cm 4cm 9cm, clip=true, width=\textwidth]{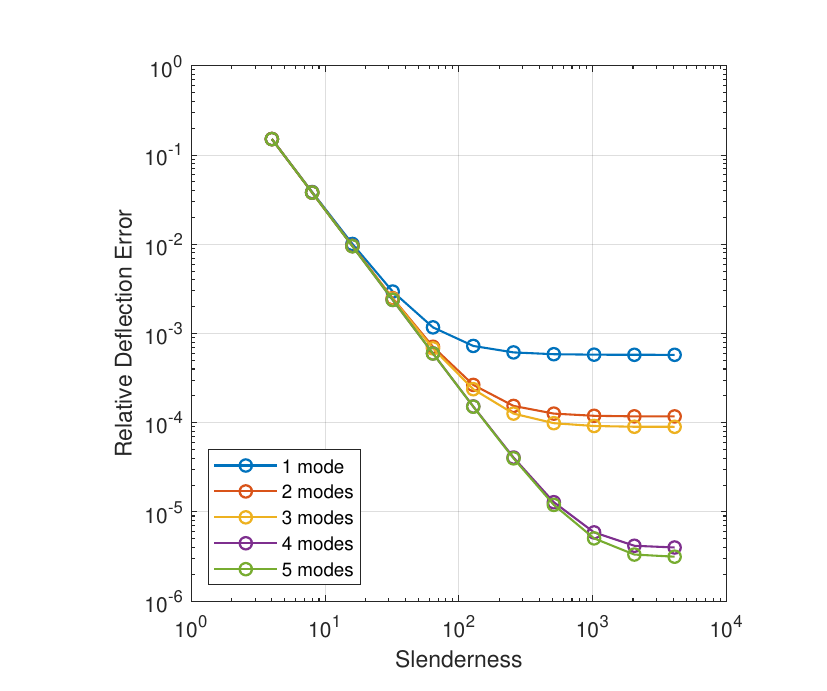}
        \caption{SS-UP}
    \end{subfigure}
    \begin{subfigure}{0.49\textwidth}
        \includegraphics[trim=3cm 9cm 4cm 9cm, clip=true, width=\textwidth]{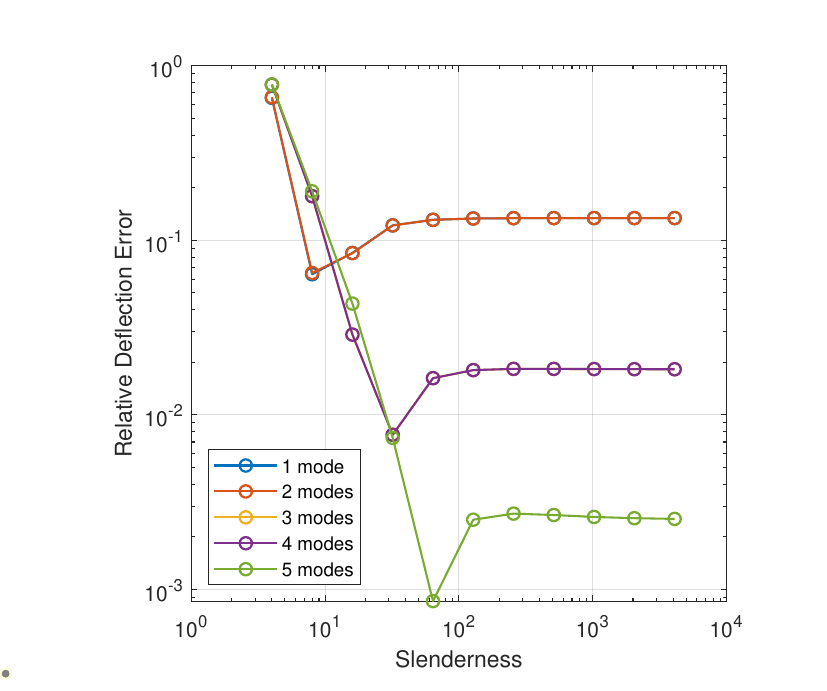}
        \caption{CC-UP}
    \end{subfigure}
    \vfill
    \begin{subfigure}{0.5\textwidth}
        \includegraphics[trim=3cm 9cm 4cm 9cm, clip=true, width=\textwidth]{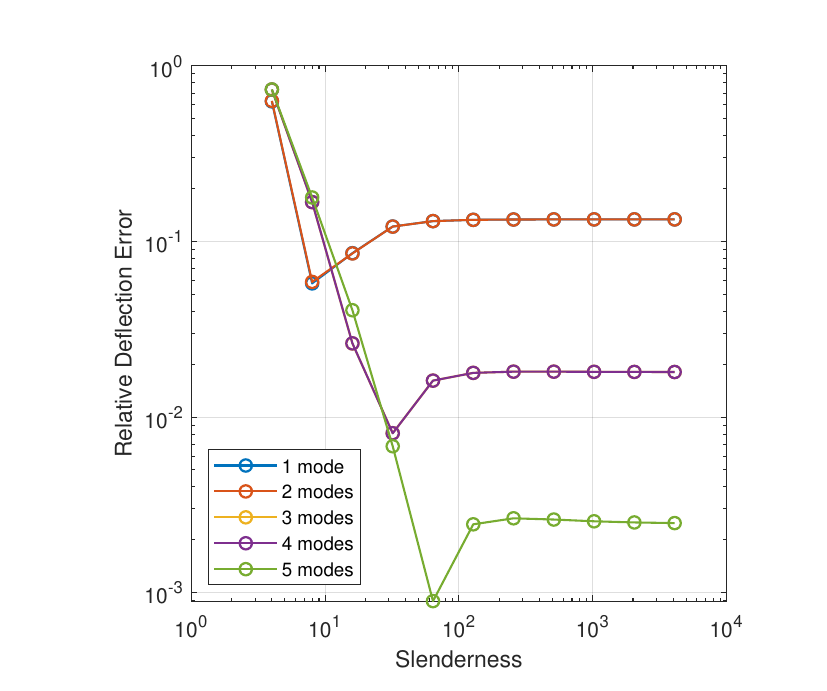}
        \caption{CC-SP}
    \end{subfigure}
    \caption{Relative deflection errors of higher-rank PGD solutions as a function of slenderness}\label{fig:CV_higher_rank_PGD}
\end{figure}

We also present results that tend to show that, regardless of the number of modes, the PGD solution is not asymptotically consistent. Figure~\ref{fig:CV_higher_rank_PGD} shows the relative deflection error as a function of slenderness for PGD solutions with multiple modes computed in a greedy way, in three different loading and boundary condition cases. The curves associated with modes 1 and 2 on the one hand, and modes 3 and 4 on the other, overlap in the CC-UP and CC-SP cases. In any case, it appears that, at least up to the use of 5 modes, the PGD solution is not consistent. However, as pointed out in Remark~\ref{rem:asymptotic_consistency}, asymptotic inconsistency does not necessarily mean that the error made by the PGD solution is large. For example, the relative deflection error is below $10^{-3}$ with a single mode in the SS-UP case. 

\subsection{Block PGD mode computation} \label{sec:new_PGD_test}

We now turn to the proposed new PGD strategy, which consists in computing the first two PGD modes as a block, i.e. simultaneously.

\subsubsection{Recovery of convergence with respect to slenderness}

\begin{figure}[h!]
    \centering
    \begin{subfigure}{0.49\textwidth}
        \includegraphics[trim=3cm 9cm 4cm 9cm, clip=true, width=\textwidth]{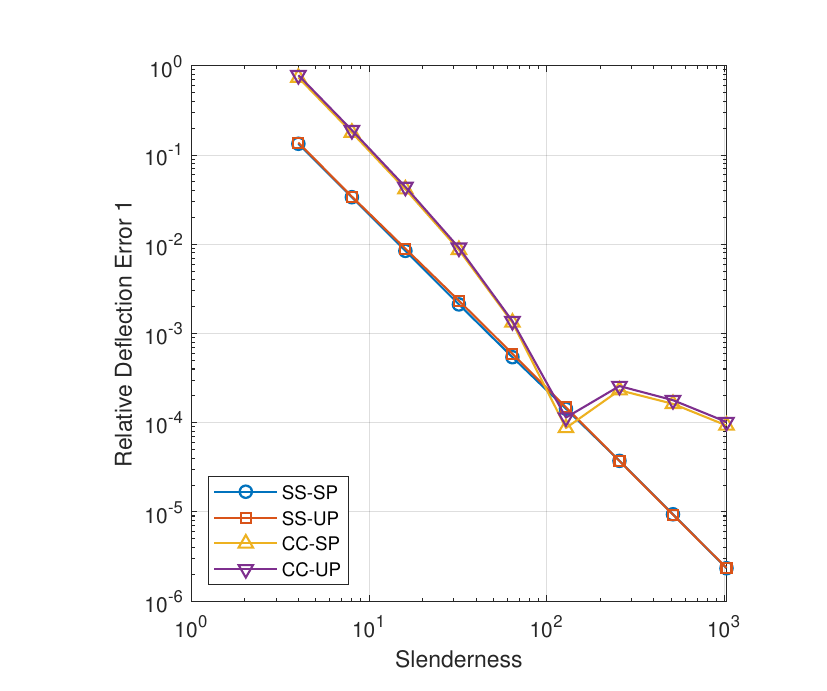}
    \end{subfigure}
    \begin{subfigure}{0.49\textwidth}
        \includegraphics[trim=3cm 9cm 4cm 9cm, clip=true, width=\textwidth]{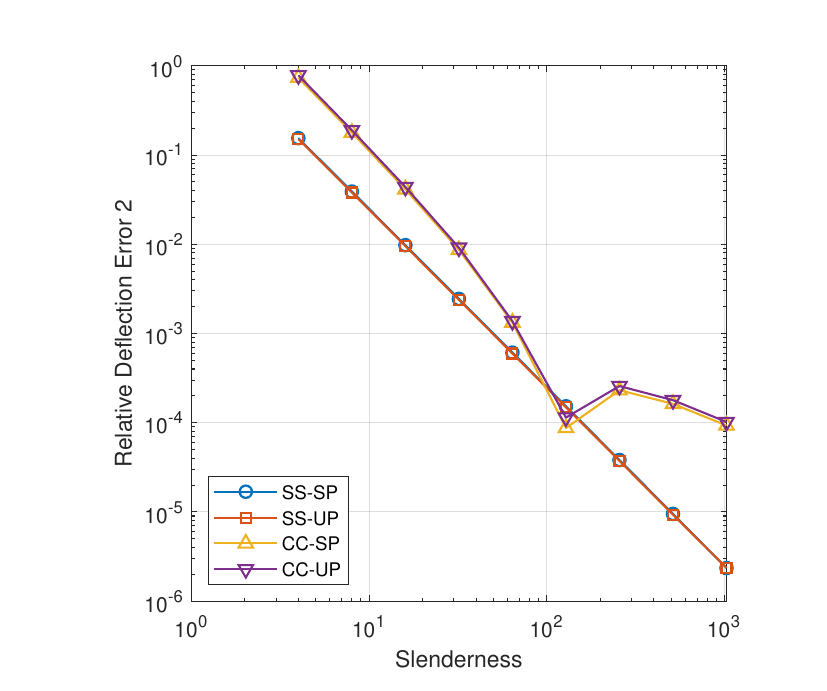}
    \end{subfigure}
    \vfill
    \begin{subfigure}{0.5\textwidth}
        \includegraphics[trim=3cm 9cm 4cm 9cm, clip=true, width=\textwidth]{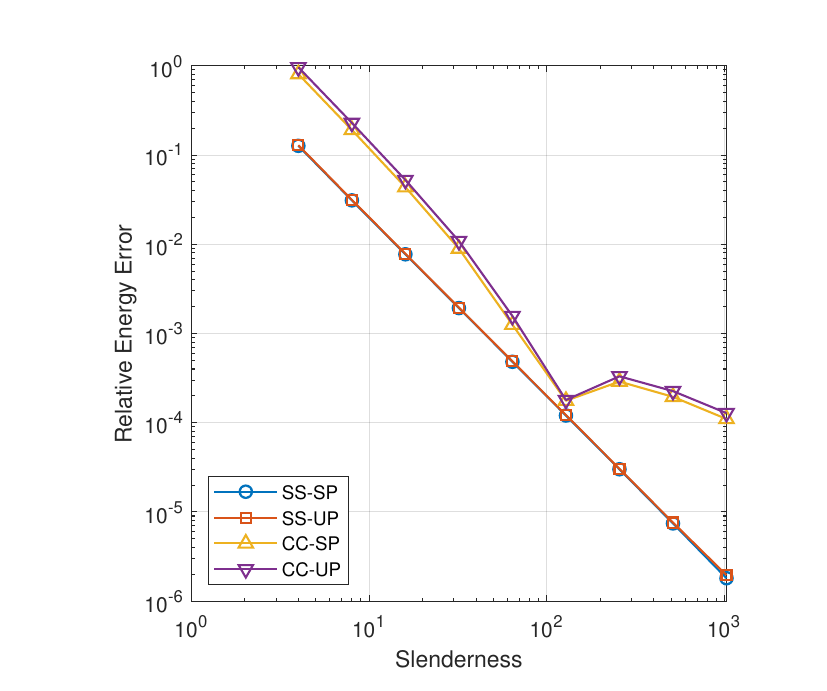}
    \end{subfigure}
    \caption{Relative deflection (top left and right) and energy (bottom) error of the first two PGD modes computed simultaneously as a function of slenderness}\label{fig:CV_new_PGD}
\end{figure}

In this section, the same convergence study is carried out as for a single PGD mode. Figure~\ref{fig:CV_new_PGD} shows the relative deflection and energy errors as a function of slenderness for two PGD modes computed simultaneously. Note that two deflection errors are shown, defined by
\begin{equation*} 
    \mathrm{Relative\ Deflection\ Error\ 1 } = \left\lvert\frac{\displaystyle \,v_3\left(\frac{L}{2}\right) - w_{KL}}{w_{KL}} \right\rvert \quad \mathrm{and} \quad  \mathrm{Relative\ Deflection\ Error\ 2} = \left\lvert\frac{\displaystyle v_3\left(\frac{L}{2}\right)+s_3(0)\,w_3\left(\frac{L}{2}\right) - w_{KL}}{w_{KL}} \right\rvert.
\end{equation*}
For the energy error, both modes ($v_3$ and $s_3\, w_3$) are taken into account in the computation of $\mathcal{E}$.

In view of these results, the simultaneous computation of this second mode enables convergence to the asymptotic solution whatever the boundary conditions and loadings. The irregularity of certain deflection error curves around a slenderness ratio of 100 reflects a change in the sign of the error. Moreover, the first two plots show that the second mode is negligible in terms of deflection, the results of the top left plot overlapping with those of the top right plot. As expected from \eqref{eq:homogenized_form}, this means that the asymptotic part of the solution is contained in the first mode. 

Finally, a decrease in the convergence rate (computed in the regime $\displaystyle 1\leq\frac{L}{t}\leq 100$) in the CC-SP and CC-UP cases is noted for slenderness larger than 100. This change in the convergence rate may be due to the presence of a boundary layer when the strip is clamped at both ends.

\subsubsection{Convergence with respect to the exact solution}

We recall that the objective of using the PGD method for plate simulation is not so much to construct new models but to approximate the 3D solution for a smaller computational cost than a direct strategy. This is why, in this section, the reference solution is no longer the Kirchhoff-Love solution, but an exact 2D solution. In practice, this reference solution is computed using finite elements on a very fine mesh. A mesh convergence test was performed. The selected mesh is sufficiently fine that the discretization error does not influence the following conclusions. We focus here on the CC-UP case. 

\begin{figure}[h!]
    \centering
    \begin{subfigure}{0.49\textwidth}
        \includegraphics[trim=3cm 9cm 4cm 9cm, clip=true, width=\textwidth]{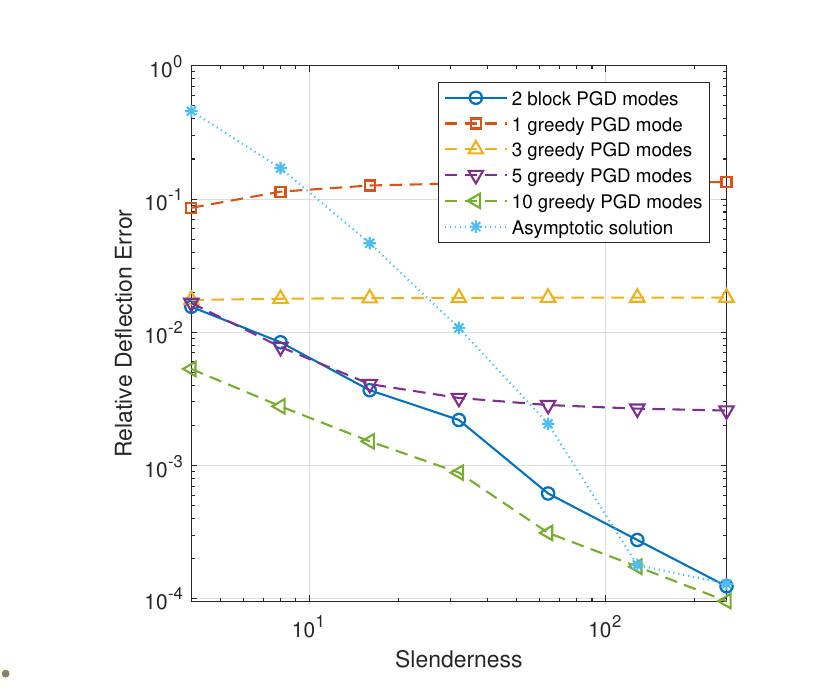}
    \end{subfigure}
    \hfill
    \begin{subfigure}{0.49\textwidth}
        \includegraphics[trim=3cm 9cm 4cm 9cm, clip=true, width=\textwidth]{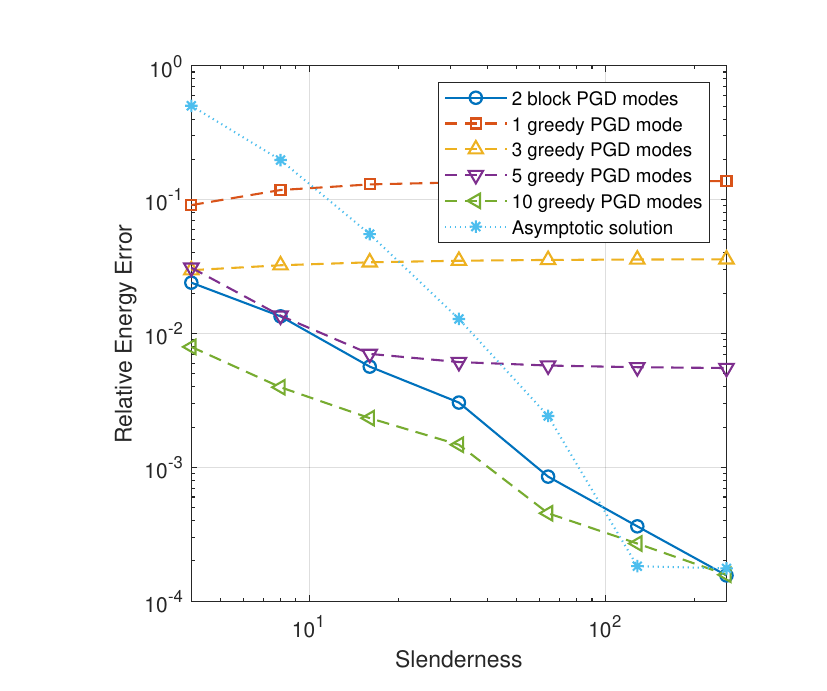}
    \end{subfigure}
    \caption{Relative deflection (left) and energy (right) error for three computational strategies: our block PGD approach simultaneoulsy computing the first two modes (solid line curves), the standard PGD approach (with different numbers of modes, dashed line curves) and the asymptotic solution (dotted line curves). The results are shown as a function of slenderness, in the CC-UP case.}\label{fig:CV_CC_UP_wrt_ref}
\end{figure}

Figure~\ref{fig:CV_CC_UP_wrt_ref} shows the deflection and energy errors as a function of slenderness for the first two modes computed simultaneously. For comparison, these errors are also shown for the standard PGD strategy with different numbers of modes and for the asymptotic solution (the 2D equivalent of \eqref{eq:homogenized_form}). The 2 block PGD modes converge to the exact solution in the limit of large slenderness, which is not the case for the standard PGD solution, at least with 5 modes or fewer. In addition, depending on the slenderness ratio, between 5 and 10 modes are required for the greedy PGD solution to be more accurate than the 2 block PGD modes. With regard to the asymptotic solution, a sign change in the error occurs around a slenderness of 100, which shows that the asymptotic regime is not reached. In any case, the 2 block PGD modes perform better than the asymptotic solution on the slenderness range of practical interest. 

\begin{figure}[h!]
    \centering
    \includegraphics[trim=3cm 9cm 4cm 9cm, clip=true, width=0.8\textwidth]{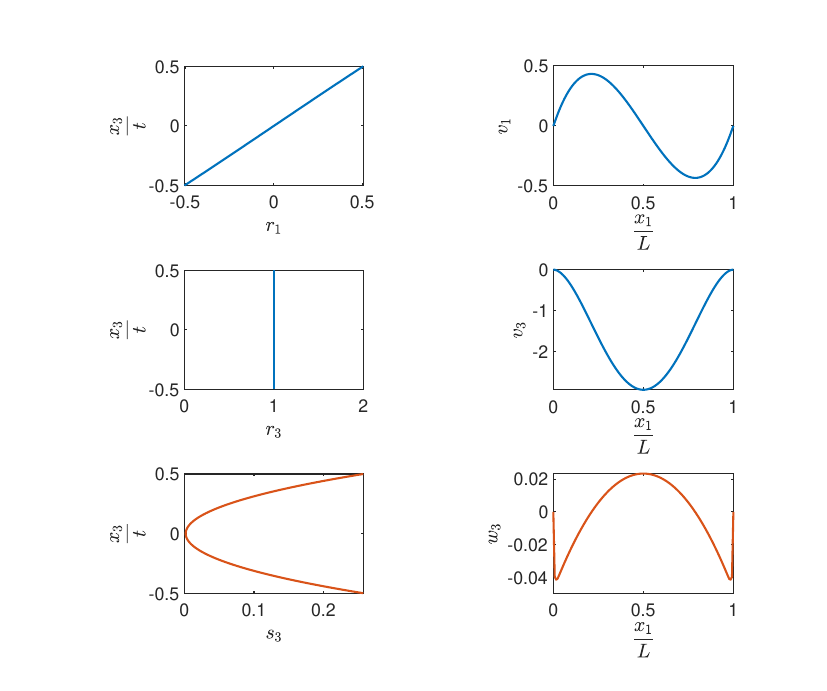}
    \caption{Components of the first two PGD modes computed simultaneously in the CC-UP case for a slenderness of 20}\label{fig:modes_CC_UP}
\end{figure}

The first two modes computed at the same time in the CC-UP case for a slenderness of 20 are shown in Figure~\ref{fig:modes_CC_UP}. The function $r_3$ is shown for visualization purposes but has not been computed. Remarkably, the first mode (in blue) exhibits a Kirchhoff-Love type kinematics: the $r_1$ function is linear and $v_1$ and $-v_{3,1}$ visually overlap. The second mode (in red) captures the boundary layer, at least in part, which is visible on the $w_3$ component. Note that the function $s_3$ is quadratic which is similar to the behavior of the corrector of the asymptotic solution (see \eqref{eq:homogenized_form}). 

\subsubsection{Addition of greedy modes}

Once the first two block PGD modes have been computed, it is of course possible to enrich the modal decomposition of the displacement if the desired accuracy is not achieved, by adding modes computed iteratively using the standard PGD algorithm. 

\begin{figure}[h!]
    \centering
    \begin{subfigure}{0.49\textwidth}
        \includegraphics[trim=3cm 9cm 4cm 9cm, clip=true, width=\textwidth]{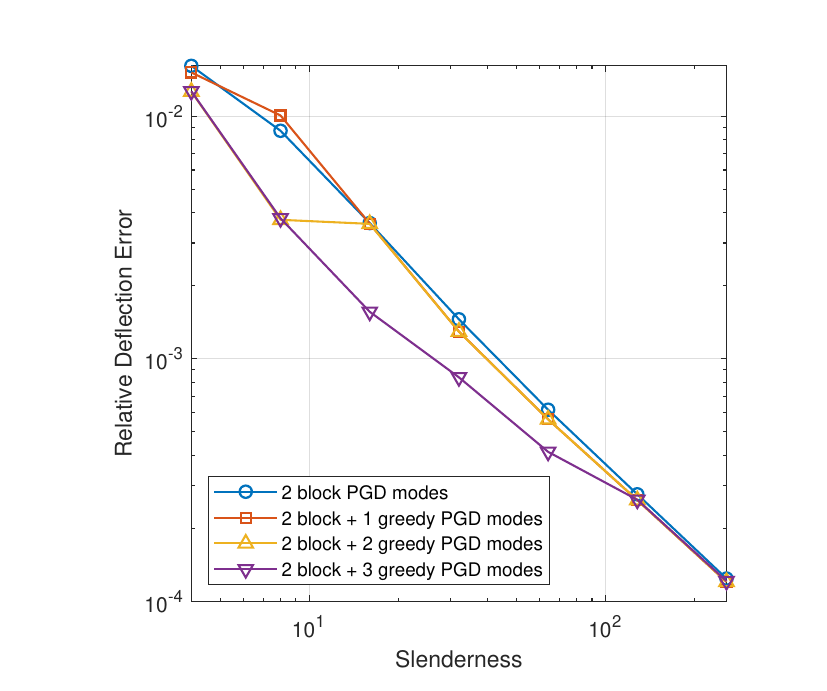}
    \end{subfigure}
    \hfill
    \begin{subfigure}{0.49\textwidth}
        \includegraphics[trim=3cm 9cm 4cm 9cm, clip=true, width=\textwidth]{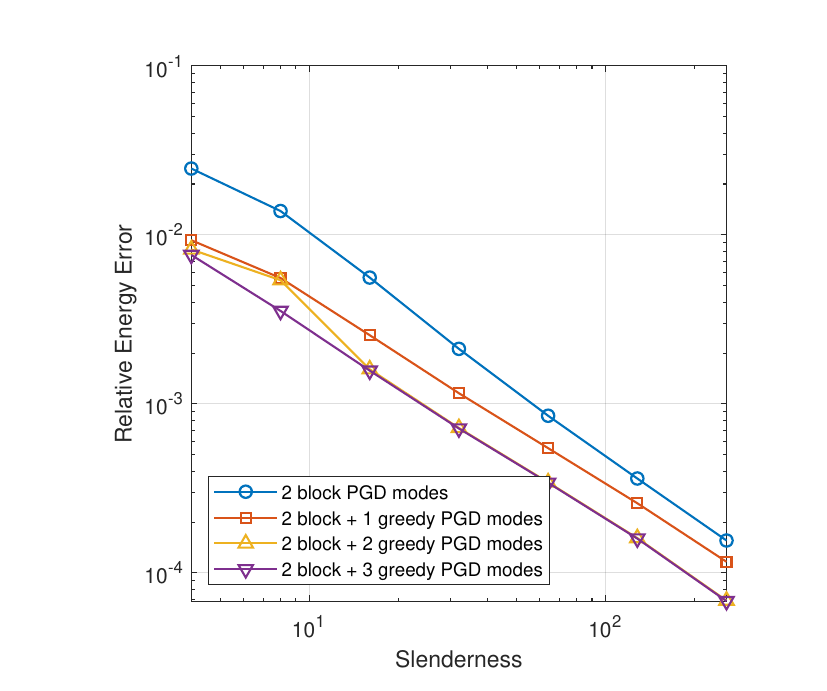}
    \end{subfigure}
    \caption{Relative deflection (left) and energy (right) error of the new PGD modal decomposition as a function of slenderness depending on the number of modes in the CC-UP case}\label{fig:CV_iterative_modes_added}
\end{figure}

Figure~\ref{fig:CV_iterative_modes_added} shows the deflection and energy errors in the CC-UP case depending on the number of modes in the PGD decomposition (these errors are computed with respect to an exact 2D solution). Depending on the number of greedy modes added, the deflection and energy estimates are improved. The relevance of modal enrichment is also dependent on the quantities of interest for the problem, the addition of modes does not affect the deflection and energy estimates in the same manner. 

\subsection{Comments on the computational cost}

In the case of bending presented here, we are aware that the new PGD strategy requires the computation of an additional function of the in-plane variable during the first iteration of the modal enrichment process. This results in an additional cost for the first PGD iteration, compared with the standard PGD method where all modes are iteratively determined. However, as highlighted in the previous section, this simultaneous computation of modes ensures that the asymptotic behavior of the exact solution is well captured, and fewer modes are needed to achieve a given accuracy. This means one can still expect to reduce the overall computational cost at the end. This point is illustrated in the CC-UP case considered before, by comparing in Table~\ref{tab:comput_cost_comp} the computational times to determine 2 block PGD modes and 5 greedy PGD modes. This comparison is largely in favor of the 2 block PGD modes. The additional cost during the first iteration with the new strategy is thus counterbalanced by the fact that a smaller number of modes is needed at the end. 

\begin{table}[h!]
    \centering
    \begin{tabular}{lccccc}
    \toprule
        & \multicolumn{5}{c}{CPU time (s)} \\
    \midrule
    Slenderness: $L/t$ & 5 & 10 & $5\times 10^1$ & $10^2$ & $2\times 10^2$ \\
    \midrule
    2 block PGD modes & 0.0547 & 0.0426 & 0.0362 & 0.0386 & 0.0491\\
    5 greedy PGD modes & 0.1434 & 0.1017 & 0.1255 & 0.1268 & 0.1604 \\
    \bottomrule
    \end{tabular}
    \caption{Comparison of CPU times for computing the first two modes simultaneously and for computing five greedy modes in the CC-UP case}
    \label{tab:comput_cost_comp}
\end{table}

\section{Conclusion}
\label{sec:conclusion}

This works compares reduced-order models of a slender strip or plate structure obtained by PGD with the solution provided by Kirchhoff-Love's theory. Using an asymptotic analysis, it is shown that, in the limit of large slenderness, the first mode of the PGD exhibits Kirchhoff-Love's kinematics but only corresponds to this model in very special cases of loading and boundary conditions. Referring to the solution given by the homogenization theory, this result can be explained by the impossibility of approximating the displacement field by a single mode in the adequate energy norm. This observation suggests a modification of the classical PGD procedure, which consists of computing several modes simultaneously. The asymptotic analysis also reveals that the PGD method is subject to locking, and we show how to deal with it by implementing a selective reduced integration technique in this context. Numerical tests show that the new PGD strategy is particularly well suited to slender structures, both in terms of accuracy and computation time. It captures the asymptotic solution as early as the first mode computation sequence and performs better than the standard PGD technique. The minimum number of block PGD modes to be computed simultaneously during the first iteration of the PGD remains an open question in the absence of material symmetries and is a perspective of this work.

\appendix

\section{Numerical implementation} \label{app1}

This appendix describes the numerical implementation associated with the results reported in Section~\ref{sec:numerical_results}. For the sake of readability, only the computation of the first approximation $\boldsymbol{u}_1$ is detailed. The determination of the modes $\boldsymbol{z}_k$ for $k\geq 2$ is standard and the presentation can easily be generalized to the 3D or laminated case. 

We recall that the first approximation $\boldsymbol{u}_1$ is sought in the form
\begin{equation*}
    \boldsymbol{u}_1(x_1,x_3) = (\boldsymbol{r}\circ \boldsymbol{v} + \boldsymbol{s}\circ \boldsymbol{w})(x_1,x_3) = 
    \begin{pmatrix}
        r_1(x_3)\,v_1(x_1) \\
        v_3(x_1) + s_3(x_3)\,w_3(x_1) 
    \end{pmatrix},
\end{equation*}
where $(v_1,v_3,w_3)$ are defined in $\mathcal{I}_1=(0,L)$ and $(r_1,s_3)$ are defined on $\displaystyle \mathcal{I}_3=\left(-\frac{t}{2},\frac{t}{2}\right)$. Using Voigt notations, we also recall that the Hooke tensor can be written in matrix form as follows:
\begin{equation*}
    \boldsymbol{C} = 
    \begin{pmatrix}
        C_{11} & C_{13} & 0 \\
        C_{13} & C_{33} & 0 \\
        0 & 0 & C_{55}
    \end{pmatrix}.
\end{equation*}

\medskip

In order to solve the PGD problem \eqref{eq:PGD2_problem} numerically, a discrete representation of the functions $v_i$, $w_3$, $r_1$ and $s_3$ is introduced. For finite element discretization in both spatial directions, we note
\begin{equation} \label{eq:FE_discr}
    v_i(x_1) = \boldsymbol{N}_1^\top(x_1)\, \boldsymbol{V}_i, \quad w_3(x_1) = \boldsymbol{N}_1^\top(x_1)\, \boldsymbol{W}_3, \quad r_1(x_3) = \boldsymbol{N}_3^\top(x_3)\, \boldsymbol{R}_1 \quad  \mathrm{and} \quad s_3(x_3)  = \boldsymbol{N}_3^\top(x_3)\,  \boldsymbol{S}_3,
\end{equation}
where $\boldsymbol{N}_1$ and $\boldsymbol{N}_3$ are the vectors of the shape functions and $(\boldsymbol{V}_i,\boldsymbol{W}_3,\boldsymbol{R}_1,\boldsymbol{S}_3)$ are the vectors of the degrees of freedom associated with each function. Alternatively, we use here a polynomial expansion in the thickness to represent the functions $r_1$ and $s_3$ \cite{vidal_proper_2013}. In this case, we also note
\begin{equation} \label{eq:PE_discr}
    r_1(x_3) = \boldsymbol{N}_3^\top(x_3)\, \boldsymbol{R}_1\quad \mathrm{and} \quad s_3(x_3)  = \boldsymbol{N}_3^\top(x_3)\,  \boldsymbol{S}_3
\end{equation}
where $\boldsymbol{N}_3(x_3)  = \begin{bmatrix} x_3^4 & x_3^3 & x_3^2 & x_3 & 1 \end{bmatrix}^\top$if we approximate $r_1$ and $s_3$ by polynomials of degree four, and $\boldsymbol{R}_1$ and $\boldsymbol{S}_3$ are the coefficients in front of each monomial function. The function $g_3$ is approximated by its interpolation on the basis of shape functions. We denote by $\boldsymbol{G}_3$ its nodal values. 

The matrices $\mathbf{K}_i$, $\mathbf{M}_i$ and $\mathbf{H}_i$, for $i\in\{1,3\}$, are next defined by
\begin{equation*}
    \mathbf{K}_i = \int_{\mathcal{I}_i} \boldsymbol{N}_i' \boldsymbol{N}_i'^\top, \quad \mathbf{M}_i = \int_{\mathcal{I}_i} \boldsymbol{N}_i \boldsymbol{N}_i^\top \quad \mathrm{and} \quad \mathbf{H}_i = \int_{\mathcal{I}_i} \boldsymbol{N}_i' \boldsymbol{N}_i^\top,
\end{equation*}
where $\boldsymbol{N}_i'$ is the vector of derivatives of the shape functions or, when using polynomial expansion in thickness, $\boldsymbol{N}_3'(x_3) = \begin{bmatrix} 4x_3^3 & 3x_3^2 & 2x_3 & 1 & 0 \end{bmatrix}^\top$.

\medskip

By expanding the two-dimensional version of equations \eqref{eq:PGD2_problem_1} and \eqref{eq:PGD2_problem_2}, separating the integrals and introducing discretization \eqref{eq:FE_discr} or \eqref{eq:PE_discr}, we obtain that $(\boldsymbol{V}_1,\boldsymbol{V}_3,\boldsymbol{W}_3)$ and $(\boldsymbol{R}_1,\boldsymbol{S}_3)$ are solutions to the following coupled matrix systems:
\begin{align}
    \begin{bmatrix}
        L_{11} & L_{12} & L_{13} \\
        L_{12} & L_{22} & L_{23} \\
        L_{13} & L_{23} & L_{33}
    \end{bmatrix}
    \begin{bmatrix}
        \boldsymbol{V}_1 \\
        \boldsymbol{V}_3 \\
        \boldsymbol{W}_3
    \end{bmatrix} &= 
    \begin{bmatrix}
        0 \\
        2\,\mathbf{M}_1\boldsymbol{G}_3 \\
        (\boldsymbol{S}_3^\top \boldsymbol{F}_3)\, \mathbf{M}_1\boldsymbol{G}_3
    \end{bmatrix}, \label{eq:algebraic_PGD_1} \\
    \begin{bmatrix} 
        P_{11} & P_{13} \\
        P_{13}   & P_{33}
    \end{bmatrix}
    \begin{bmatrix}
        \boldsymbol{R}_1 \\
        \boldsymbol{S}_3
    \end{bmatrix} &= 
    \begin{bmatrix}
        0 \\
        (\boldsymbol{W}_3^\top\mathbf{M}_1 \boldsymbol{G}_3)\, \boldsymbol{F}_3
    \end{bmatrix} - 
    \begin{bmatrix}
        C_{55}\,(\boldsymbol{V}_3^\top \mathbf{H}_1\boldsymbol{V}_1)\,\mathbf{H}_3\,\boldsymbol{R}_3 \\
        C_{55}\,(\boldsymbol{W}_3^\top \mathbf{K}_1\mathbf{V}_3)\,\mathbf{M}_3\,\boldsymbol{R}_3
    \end{bmatrix}, \label{eq:algebraic_PGD_2}
\end{align}
where $\boldsymbol{F}_3$ is a vector such that $\boldsymbol{S}_3^\top \boldsymbol{F}_3 = s_3^+ + s_3^-$, $\boldsymbol{R}_3$ is the discrete representation of the constant function equal to 1 on $\mathcal{I}_3$ and
\begin{equation} \label{eq:PGD_matrices}
    \begin{aligned}
    L_{11} &= C_{11}\,(\boldsymbol{R}_1^\top \mathbf{M}_3\boldsymbol{R}_1)\, \mathbf{K}_1 + C_{55}\,(\boldsymbol{R}_1^\top \mathbf{K}_3\boldsymbol{R}_1)\,\mathbf{M}_1, \\
    L_{22} &= C_{55}\,t\,\mathbf{K}_1, \\
    L_{33} &= C_{33}\,(\boldsymbol{S}_3^\top \mathbf{K}_3\boldsymbol{S}_3)\,\mathbf{M}_1 + C_{55}\,(\boldsymbol{S}_3^\top \mathbf{M}_3\boldsymbol{S}_3)\,\mathbf{K}_1, \\
    L_{12} &= C_{55}\,(\boldsymbol{R}_1^\top \mathbf{H}_3\boldsymbol{R}_3)\,\mathbf{H}_1^\top, \\
    L_{13} &= C_{13}\,(\boldsymbol{S}_3^\top \mathbf{H}_3\boldsymbol{R}_1)\,\mathbf{H}_1 + C_{55}\,(\boldsymbol{R}_1^\top \mathbf{H}_3\boldsymbol{S}_3)\,\mathbf{H}_1^\top, \\
    L_{23} &= C_{55}\,(\boldsymbol{S}_3^\top \mathbf{M}_3\boldsymbol{R}_3)\,\mathbf{K}_1, \\
    P_{11} &= C_{11}\,(\boldsymbol{V}_1^\top \mathbf{K}_1\boldsymbol{V}_1)\,\mathbf{M}_3 + C_{55}\,(\boldsymbol{V}_1^\top \mathbf{M}_1\boldsymbol{V}_1)\,\mathbf{K}_3, \\
    P_{33} &= C_{33}\,(\boldsymbol{W}_3^\top \mathbf{M}_1\boldsymbol{W}_3)\,\mathbf{K}_3 + C_{55}\,(\boldsymbol{W}_3^\top \mathbf{K}_1\boldsymbol{W}_3)\,\mathbf{M}_3, \\
    P_{13} &= C_{13}\,(\boldsymbol{V}_1^\top \mathbf{H}_1\boldsymbol{W}_3)\,\mathbf{H}_3^\top + C_{55}\,(\boldsymbol{W}_3^\top \mathbf{H}_1\boldsymbol{V}_1)\,\mathbf{H}_3.
\end{aligned}
\end{equation}
The coupled systems \eqref{eq:algebraic_PGD_1}-\eqref{eq:algebraic_PGD_2} are solved using the fixed-point algorithm described in the continuous case at the end of Section~\ref{sec:new_PGD_strategy}. The vectors $\boldsymbol{R_1}$ and $\boldsymbol{S}_3$ are also normalized at the end of each iteration of the fixed point algorithm. 

To prevent shear locking, we resort to selective integration to under-integrate parts of the stiffness matrices associated with transverse shear, i.e. to integrate with fewer integration points than necessary. The resulting modifications to \eqref{eq:PGD_matrices} are given below:
\begin{equation*}
\begin{aligned}
     L_{11} &= C_{11}\,(\boldsymbol{R}_1^\top \mathbf{M}_3\boldsymbol{R}_1)\,\mathbf{K}_1 + C_{55}\,(\boldsymbol{R}_1^\top \mathbf{K}_3\boldsymbol{R}_1)\,\mathbf{M}_1^{RI} \\   
     P_{11} &= C_{11}\,(\boldsymbol{V}_1^\top \mathbf{K}_1\boldsymbol{V}_1)\,\mathbf{M}_3 + C_{55}\,(\boldsymbol{V}_1^\top \mathbf{M}_1^{RI}\boldsymbol{V}_1)\, \mathbf{K}_3 
\end{aligned}
\end{equation*}
where $\mathbf{M}_1^{RI}$ is obtained using a reduced integration rule. For quadratic Lagrange elements, 2 Gauss points per element are used instead of 3. 

%% For citations use: 
%%       \cite{<label>} ==> [1]

%%

%% If you have bib database file and want bibtex to generate the
%% bibitems, please use
%%
%%  \bibliographystyle{elsarticle-num} 
%%  \bibliography{<your bibdatabase>}

\bibliographystyle{elsarticle-num} 
\bibliography{references}

%% else use the following coding to input the bibitems directly in the
%% TeX file.

%% Refer following link for more details about bibliography and citations.
%% https://en.wikibooks.org/wiki/LaTeX/Bibliography_Management

%%\begin{thebibliography}{00}

%% For numbered reference style
%% \bibitem{label}
%% Text of bibliographic item

%%\bibitem{lamport94}
  %%Leslie Lamport,
  %%\textit{\LaTeX: a document preparation system},
  %%Addison Wesley, Massachusetts,
  %%2nd edition,
  %%1994.

%%\end{thebibliography}

\end{document}